\documentclass[11pt]{amsart}
\usepackage[margin=1in]{geometry}
\usepackage[T1]{fontenc}
\usepackage{lmodern,amsfonts,amsmath,amstext,amsbsy,amssymb,amsopn,amsthm,upref,eucal,mathptmx,mathtools,url,thmtools
}
\usepackage[foot]{amsaddr}
\IfFileExists{dsfont.sty}{\usepackage{dsfont}}{}
\providecommand{\mathds}[1]{\mathbf{#1}}
\usepackage{enumitem,xcolor}
\usepackage[colorlinks,hypertexnames=false]{hyperref}
\usepackage{cleveref}
\newtheorem{thm}{Theorem}[section]
\newtheorem{prop}[thm]{Proposition}
\newtheorem{defn}[thm]{Definition}
\newtheorem{lem}[thm]{Lemma}
\newtheorem{oss}[thm]{Remark}
\newtheorem{ex}[thm]{Example}
\newtheorem{cor}[thm]{Corollary}
\numberwithin{equation}{section}
\newtheorem{bigthm}{Theorem}

\definecolor{auditgreen}{RGB}{0,115,75}
\definecolor{revisionblue}{RGB}{0,70,190}
\colorlet{REVISIONBLUE}{revisionblue}

\newcommand{\dd}{\mathrm{d}}
\newcommand{\fF}{\mathfrak{F}}

\newcommand\cA{{\mathcal A}}
\newcommand\cB{{\mathcal B}}
\newcommand\cC{{\mathcal C}}

\newcommand\cE{{\mathcal E}}

\newcommand\cH{{\mathcal H}}

\newcommand\cL{{\mathcal L}}

\newcommand\cN{{\mathcal N}}

\newcommand\cT{{\mathcal T}}
\newcommand\cU{{\mathcal U}}

\newcommand\bC{{\mathbb C}}

\newcommand\bN{{\mathbb N}}

\newcommand\bP{{\mathbb P}}

\newcommand\bR{{\mathbb R}}
\newcommand\bT{{\mathbb T}}

\newcommand\bZ{{\mathbb Z}}

\newcommand\ve{\varepsilon}

\newcommand\frh{{\mathfrak h}}

\newcommand{\frm}{\mathfrak{m}}

\newcommand\fC{{\mathbf{C}}}
\newcommand\fn{{\mathfrak{n}}}

\newcommand{\supp}{\operatorname{supp}}

\newcommand\Id{\text{Id}}
\begin{document}
\title[Dolgopyat estimates for partially hyperbolic systems]{{Dolgopyat estimate for partially hyperbolic systems: limit theorems and exponential mixing}}
\author{Roberto Castorrini \textsuperscript{(1)}}
\address{{\tiny (1)} Dipartimento di Economia, Ingegneria, Società e Impresa (DEIM), Università della Tuscia, 01100, Viterbo, Italy.}
\email{roberto.castorrini@unitus.it}
\author{Kasun Fernando \textsuperscript{(2)}}
\address{{\tiny (2)} Brunell University of London, Kingston Ln, Uxbridge UB8 3PH, UK.}
\email{{\tt kasun.fernando@brunel.ac.uk }}
\thanks{This paper started when RC was working at Scuola Normale Superiore di Pisa, and he was partially supported by the research project ‘Dynamics and Information Research Institute—Quantum Information, Quantum Technologies’ within the agreement between UniCredit Bank and Scuola Normale Superiore and by the research project PRIN 2017S35EHN of the MIUR. The authors acknowledge the UMI Group “DinAmicI” (\url{www.dinamici.org}) and the INdAM group GNFM. 
Finally, it is a pleasure to thank the three Cs, Claudio Bonanno, Carlo Collari and Carlangelo Liverani, for helpful discussions. \\
} 

\begin{abstract}{We prove a Dolgopyat-type estimate for the twisted transfer operators of a class of two-dimensional partially hyperbolic endomorphisms of the torus, not necessarily skew-products. The twist satisfies a quantitative transversality condition on the cotangent dynamics of the associated three-dimensional extension which, in the skew-product case, is related to Tsujii's notions of transversality and partial captivity conditions. We use the main estimate to prove higher order expansions in the central and mixing local limit theorems for the two-dimensional base, including the case of several ergodic absolutely continuous invariant measures, and we prove exponential decay of correlations for the associated three-dimensional circle extensions, with a partially hyperbolic base.}
\end{abstract}
\maketitle

\section{Introduction}
Quantitative mixing estimates are an important tool in the study of statistical limit laws, including central and local limit theorems, Berry--Esseen estimates, and
large deviation principles. For systems with a neutral direction, proving exponential mixing often requires a cancellation mechanism, first exploited by Dolgopyat for Anosov flows \cite{Dima98} and compact group extensions of hyperbolic systems \cite{Dima02}.
For a circle extension $(x,z)\mapsto(f(x),z+\tau(x))$ over an expanding map, Fourier decomposition in the fibre leads to the twisted transfer operators $\cL_\nu=\cL_f(e^{i\nu\tau}\,\cdot\,)$. 
A Dolgopyat-type inequality gives quantitative bounds on
the twisted transfer operators at large frequencies.
Together with suitable spectral information at bounded
frequencies, these bounds yield exponential mixing and
higher-order statistical expansions, as we show below.
For results on mixing and related spectral
estimates for extensions of expanding maps see
\cite{Dima02,Gou09,F11,NTW,ChHu,But24};
for suspension semiflows, see \cite{Tsu08}.
The cancellation method underlying Dolgopyat-type
estimates originates in \cite{Dima98} and it was used even in the presence of discontinuities \cite{BuEs}. However, if the system does not have a skew-product structure, the mixing problem becomes much harder to treat, and even proving existence of invariant measures is highly non-trivial. 

The main technical result of this paper is the Dolgopyat estimate (Theorem~\ref{thm:DI}) for a broad class of two-dimensional partially hyperbolic system for which existence of finitely many invariant mixing densities is known. We consider the strongly dominated vertical partially hyperbolic maps (SVPHs) of \cite{CaLi}, restricted to a class of virtually expanding maps in the sense of Tsujii \cite{Tsu1}. The latter restriction is made to work with the classical Sobolev spaces $H^s(\bT^2)$. Extending the argument to the anisotropic spaces of \cite{CaLi} is possible, but it requires additional operator and localisation estimates. The main difficulty faced by the paper is that inverse branches need not contract the central direction, so the standard expanding-base arguments do not apply directly. Our result apply in particular to the extension 
\[
\tilde F(x,\theta,z)=(F(x,\theta),z+\tau(x,\theta))\pmod{\bZ^3},
\]
where $F$ is partially hyperbolic, under a quantitative cotangent transversality condition analogous to the partial captivity conditions in \cite{F11,NTW}.\\
We give two important applications of Theorem~\ref{thm:DI}. Theorem~\ref{teoA} establishes Edgeworth expansions and higher order mixing global and local limit expansions for the Birkhoff sums of $\tau$ over $F$. Several ergodic acips are allowed and the leading distribution is a mixture of centred Gaussians, with weights determined by the initial density. The spectral gap and perturbation near zero give the Gaussian terms, while the Dolgopyat estimate controls the large real frequencies needed for expansions of arbitrary order. This application uses the spectral methods developed in \cite{HeHe,FeLi,FePe}. Theorem~\ref{teoB} establishes exponential mixing of $\tilde F$ for H\"older observables.

For skew products $F(x,\theta)=(f(x),\theta+\omega(x))$ with $\tau=\tau(x)$, the extension is a $\bT^2$-extension of $f$ with cocycle $(\omega,\tau)$. In that case exponential mixing is already known from \cite{ChHu} when $a\omega+b\tau$ is not cohomologous to a constant for every $(a,b)\ne(0,0)$. Lemma~\ref{lem:skew} identifies \textbf{(B1)} with this condition (see also Remark~\ref{rmk:skew-sharp} for its necessity for DI at all large real frequencies). Our estimate also applies when $F$ is not a skew product, or when $\tau$ depends on the central coordinate, provided the stated hypotheses hold.

The proof combines Lasota--Yorke inequalities with an $L^2$ estimate at times of order $\log|\nu|$. We localise in space at scale $|\nu|^{-1/2}$ and in frequency. Condition \textbf{(B1)} controls frequencies reached by many inverse branches and non-stationary phase controls the remaining frequencies.

Section~\ref{sec:setting} gives the assumptions, the technical theorem, and its two applications. Section~\ref{sec:transfer} develops the spectral tools. Section~\ref{sec:proof-main} proves the Dolgopyat estimate, followed by the proofs of Theorems~\ref{teoA} and~\ref{teoB} in Sections~\ref{sec:proof-thmA} and~\ref{sec:doc}.  Appendix~\ref{app:B1} discuss the main assumption in the skew-roduct case and provide a simple example to which our results apply and, finally, Appendix~\ref{sec:previous} recalls the needed estimates from \cite{CaLi}.

\section{Setting and main results}\label{sec:setting}

\subsection{Assumptions on the map}
We first introduce the systems $F\in \cC^{\infty}(\bT^2,\bT^2)$ studied in \cite{CaLi} and called {\it strongly dominated vertical partially hyperbolic systems}, or SVPHs for short.\footnote{Here we require the maps to be smooth even though in \cite{CaLi} $r \geq 4$ suffices.} For the convenience of the reader we recall the main assumptions. A set of sufficient conditions for checking whether a system is an SVPH are provided in \cite[Lemma 2.2]{CaLi}. We refer to \cite[Section 2.2]{CaLi} for a detailed discussion, including examples, of assumptions below, and we refer to the Appendix \ref{sec:previous} for a collection of technical results about SVPHs.\\

Let $F:\mathbb{T}^2\to \mathbb{T}^2$ be a surjective $\mathcal{C}^{\infty}$ local diffeomorphism of the form
\begin{equation}
\label{the map F}
{F(x,\theta)=(f(x,\theta),\theta +\omega(x,\theta)) \mod \bZ^2,\qquad \omega\in\cC^\infty(\bT^2,\bR).}
\end{equation}
We call $F$ an SVPH if $F$ satisfies the following properties \textbf{(A1)}--\textbf{(A5)}.\\

\textbf{(A1)} For all $p\in\bT^2$ we have $\det(D_pF)>0$.\\

\textbf{(A2)} There exist $\chi_u\in (0,1)$, $\chi_c\in (0,1]$ and $0<\lambda^{-}_c< 1< \lambda_c^+ < \lambda\le\Lambda$ such that, setting
\begin{equation}\label{eq:cone-def}
\begin{split}
&\fC_u:= \lbrace (\xi, \eta)\in \cT_z\mathbb{T}^2  :  |\eta|\le \chi_u|\xi | \rbrace \\
&\fC_c:= \lbrace (\xi, \eta)\in  \cT_z\mathbb{T}^2 :  |\xi|\le \chi_c |\eta | \rbrace,
\end{split}
\end{equation}
defining
\begin{equation}\label{def of lamb^+ mu^+}
\begin{split}
&\lambda^-_n(z):= \inf_{v\in \bR^2\setminus \fC_c }\mathfrak{m}_{F^{n}}(z, v) \qquad \qquad\qquad\; \lambda^+_n(z):= \sup_{v\in \bR^2\setminus \fC_c }\mathfrak{m}_{F^{n}}(z, v),\\
&\lambda_{c,n}^-(z):= \inf_{v\in \fC_c\setminus{\lbrace 0 \rbrace}}\mathfrak{m}^*_{F^{n}}(z, v) \quad\quad\qquad\qquad\lambda_{c,n}^+(z)
:= \sup_{v\in \fC_c\setminus{\lbrace 0 \rbrace}}\mathfrak{m}^*_{F^{n}}(z, v),
\end{split}
\end{equation}
where $\mathfrak{m}_F(z, v)=\frac{\|D_zF v\|}{\|v\|}$ for $v\in \cT_z\bT^2$, and $\mathfrak{m}^*_F(z, v)=\frac{\|(D_zF)^{-1} v\|}{\|v\|}$ for $v\in \cT_{F(z)}\bT^2$,
and letting $\lambda^-_n=\inf_z\lambda^-_n(z)$ and $\lambda^+_n=\sup_z\lambda^+_n(z)$ we assume the following:\\
There exist constants $C_\star\geq 1$ and $\iota_\star \in (0,1)$ such that, for all $z\in\bT^2$ and $n\in\bN$,\footnote{$A\Subset B$ means $\overline{A}\subset \mbox{int}(B)\cup \lbrace 0 \rbrace$.}
\begin{equation}
  \begin{aligned}
  D_zF\fC_u &\subset \{(\xi,\eta): |\eta|\le \iota_\star\chi_u|\xi|\}\Subset \fC_u; \\ (D_zF)^{-1}\fC_c &\subset \{(\xi,\eta): |\xi|\le \iota_\star\chi_c|\eta|\}\,, \label{invariance of cone}
  \end{aligned}
\end{equation}
\begin{equation}
  \begin{aligned}
  C_\star^{-1}({\lambda_c^-})^n&\leq \lambda_{c,n}^-(z)\leq \lambda_{c,n}^+(z)\leq C_\star({\lambda_{c}^+})^n \;;\\ C_\star^{-1}\lambda^n&\leq\lambda^-_n(z)\leq  \lambda^+_n(z)\leq C_\star \Lambda^n \,.\label{partial hyperbolicity 2}
  \end{aligned}
\end{equation}

\vspace{15pt}
\textbf{(A3)} Let $\Upsilon$ be the family of closed curves $\gamma \in \cC^\infty(\bT, \bT^2)$ such that\footnote{As usual we consider equivalent two curves that differ only by a $\cC^\infty$ reparametrization. In the following, we will mostly use curves that are parametrized by vertical length.}
\begin{itemize}
\item[c0)] $\gamma'\neq 0$,
\item[c1)] $\gamma$ has homotopy class $(0,1)$,
\item[c2)] $\gamma'(t)\in \fC_c$, for each $t\in \bT$.
\end{itemize}
We assume that for each $\gamma\in\Upsilon$ there exist $\{\gamma_i\}_{i=1}^N\subset \Upsilon$ such that $F^{-1}(\gamma)=\bigcup_{i=1}^N\gamma_i$.\\

\textbf{(A4)} We assume that $F$ satisfies the pinching condition\footnote{In \cite{CaLi} there is a much stronger pinching condition due to the Banach norms used. Here we just need to assume a dominated splitting.}
\begin{equation}\label{pinching condition}
\max\{\lambda_c^+,\,(\lambda_c^-)^{-1}\} <\lambda.
\end{equation}

\vspace{15pt}
\textbf{(A5)} With the notation \eqref{the map F} we require, for each $p\in\bT^2$,
\[
\partial_xf(p)>\max\left\{ 2(1+\|\partial_x\omega\|_\infty),|\partial_\theta f(p)|\right\}.
\]

From now on we set
\begin{equation}\label{eq:mumax}
{\lambda_c}:=\max\{ \lambda_c^+,1/\lambda^-_c\}\ge  1,
\end{equation}
so that \eqref{pinching condition} reads $\lambda_c<\lambda$.\\

Next, the spectral gap for the transfer operator associated to $F$ allows for nontrivial peripheral eigenvalues,
corresponding to finite cyclic decompositions of the ergodic components. To state our statistical applications without
passing to iterates, we impose the following assumption besides \textbf{(A1)}--\textbf{(A5)} .

\medskip
\noindent\textbf{(A6)}
Every ergodic absolutely continuous invariant probability
measure of $F$ is weakly mixing.\\

Finally, note that since $F$ is a local diffeomorphism of the compact torus,
so we can fix constants $C_\infty\ge1$ and $\Theta\ge1$ such that
\begin{equation}\label{eq:Linfty-growth}
\|\cL_F^n1\|_\infty\le C_\infty\Theta^n
\qquad\text{for all }n\in\bN.
\end{equation}

\begin{ex}
Notable examples of SVPH are partially expanding systems with neutral or close to neutral central direction which includes:
\begin{itemize}[leftmargin=*]
\item Skew torus extensions over smooth expanding map (and small perturbations of it): $F(x,\theta)=(f(x), \theta+\omega(x))$ on $\bT^2$ where $f$ is expanding.
\item Fast-slow partially expanding systems: $F_\ve(x,\theta)=(f(x,\theta), \theta+ \ve \omega(x,\theta))$ on $\bT^2$ where $\inf_{(x,\theta)\in \bT^2}\partial_x f(x,\theta)>1$ and $\ve>0$ is small, under the further assumptions of \cite[Section 10]{CaLi}.
\end{itemize}
One can check that the first example satisfies the assumptions in \cite[Lemma 2.2]{CaLi} and hence is an SVPH, provided the expansion is large enough.{\footnote{For a general cocycle, passing to an iterate alone does not guarantee the coordinate inequalities in \textbf{(A5)}; these must still be checked.}} For a detailed discussion of the second, see \cite[Section 10]{CaLi} and \cite{DeLi2, Liv}.
\end{ex}

The main statistical properties of SVPHs (existence of invariant measures, regularity, decay of correlations) are proved in \cite{CaLi} under a transversality property of the unstable cones, which is generic \cite{Tsu}. We recall it, as it will appear in the case of skew products (Lemma \ref{lem:skew}).

\begin{defn}
\label{transversality}
Given $n\in \mathbb{N}$, $y\in \mathbb{T}^2$ and $z_1,z_2 \in F^{-n}(y)$, we say that $z_1$ is \textit{transversal} to $z_2$ (at time $n$) if $D_{z_1} F^n \fC_u \cap D_{z_2} F^n \fC_u=\lbrace 0 \rbrace$, and we write $z_1 \pitchfork_n z_2$. We set
\[
\mathcal{N}_{F}(n):=\sup_{y\in \mathbb{T}^2}\ \sup_{z_1\in F^{-n}(y)}\ \sum_{\substack{z_2\in F^{-n}(y)\\ z_2 \not{\pitchfork}_n z_1 }}  |\det D_{z_2}F^{n}|^{-1},
\]
and we say that $F$ satisfies the \emph{transversality condition} if
\begin{equation}\label{main-hp}\tag{T}
\text{there exist $C>0$, $\gamma_0>0$ such that}\quad \cN_F(n)\le Ce^{-\gamma_0 n}\quad\text{for all }n\in\bN .
\end{equation}
\end{defn}

\vspace{10pt}
\begin{oss}\label{rmk:A8}
The constants in \eqref{eq:Linfty-growth}  always exist. The relevant restriction is the compatibility condition \textbf{(P)} below, which requires an admissible growth rate $\Theta$
sufficiently close to $1$ relative to the spectral and transversality rates.
By \eqref{C1norm}, \eqref{eq:Linfty-growth} holds for every SVPH with $\Theta=\lambda_c$. For skew products $F(x,\theta)=(f(x),\theta+\omega(x))$ over expanding maps it holds with $\Theta=1$, since then $\sum_{y\in F^{-n}(x,\theta)}|\det D_yF^n|^{-1}=\sum_{x'\in f^{-n}(x)}|(f^n)'(x')|^{-1}$, which is bounded uniformly in $n$. In our results, $\Theta$ enters only through the condition \textbf{(P)} below, which requires $\Theta$ to be close to $1$ in terms of the other rates.
\end{oss}

\subsection{Virtually expanding maps}

We restrict to the class of systems introduced in \cite{Tsu1}. This allows us to work directly with the classical Sobolev spaces $H^s(\bT^2)$, rather than with the anisotropic Banach spaces used in \cite{CaLi}, although the arguments can be adapted to this case.

We consider the push-forward action of $F$ on the cotangent space $\cT^*\bT^2=\bT^2\times \bR^2$ given by
\[
F^{\dagger}(p, \eta)=(F(p), [(D_pF)^*]^{-1}\eta),
\]
where $^*$ denotes the transpose of a matrix. Note that $F^\dagger$ is a local diffeomorphism and a finite covering of $\cT^*\bT^2$, and $(F^n)^\dagger=(F^\dagger)^n$.

For $s \ge 0$ we consider the quantity
\begin{equation}\label{Bs}
B_s(F):=\sup_{q\in\bT^2,\ \xi\in\bR^2\setminus\{0\}}\ \sum_{F^{\dagger}(p, \eta)=(q, \xi)}\left( \frac{\|\xi\|}{\|\eta\|}\right)^{s}\frac{1}{|\det D_pF|}.
\end{equation}

{We assume that there is $s>0$ such that $B_s(F^n)$ decays exponentially fast enough to compensate for the possible expansion of inverse branches in the central direction:}
\begin{equation}\label{virtually}
\lambda_c^s \inf _{n \geq 1}\left(B_s\left(F^{n}\right)\right)^{1 / n}=\lambda_c^s\lim _{n \rightarrow \infty}\left(B_s\left(F^{n}\right)\right)^{1 / n}<1,
\end{equation}
where we have used submultiplicativity. Following the naming convention in \cite{Tsu1}, we refer to an SVPH $F$ which satisfies \eqref{virtually} as $s$-\textit{Virtually Expanding} ($s$-VE).\footnote{In \cite{Tsu1}, $F$ is called $s$-virtually expanding if $\lim_{n \rightarrow \infty}B_{2s}(F^{n})^{1 / n}<1$. Condition \eqref{virtually} is stronger, see Lemma \ref{lem:VE-Tsujii}.} When the role of $s$ is not important, we will drop the prefix and call such maps virtually expanding, or simply VE.

\begin{ex}
In \cite{Tsu1} it is shown that suitable skew products $F(x,\theta)=(f(x), \theta+\omega(x))$ over expanding maps of the circle, with sufficiently large expansion, are virtually expanding. We refer to \cite{Tsu1} for the precise statement and for further examples.
\end{ex}

\subsection{Assumptions on the observable \texorpdfstring{$\tau$}{Lg}}

We will consider observables $\tau \in \cC^\infty(\bT^2,\bR)$.\footnote{Finitely many derivatives suffice, but for simplicity we do not track the number.} For $q\in\bT^2$ and $n\in\bN$ we denote by $\cH^n(q)$ the set of inverse branches of $F^n$ defined in a neighbourhood of $q$, and set $J_\frh(q):=|\det D_{\frh(q)}F^n|^{-1}$ for $\frh\in\cH^n(q)$. Recall $\tau_n = \sum_{k=0}^{n-1} \tau \circ F^k$. For a differentiable function $\varphi$, we denote by
$d\varphi_q$ its differential at $q$. Thus, for
$u\in T_q\bT^2$,
\[
d(\tau_n\circ\frh)_q(u)
=\nabla(\tau_n\circ\frh)(q)\cdot u
=d(\tau_n)_{\frh(q)}(D_q\frh\,u).
\]
We assume the follwing\\
\vspace{15pt}
\textbf{(B1)} There exist $C>0$, $\gamma_1>0$ and $\kappa\in(0,\log\lambda)$ such that, for all $n\in\bN$, $q\in\bT^2$ and $\zeta\in(\bR^2)^*$,
\begin{equation}\label{eq:B1star}
\sum_{\frh\in\cH^n(q)\,:\ \exists\, u\in D_{\frh(q)}F^n\fC_u,\ \|u\|=1,\ |\zeta(u)-d(\tau_n\circ\frh)_q(u)|\le e^{\kappa n}\lambda^{-n}}J_\frh(q)\ \le\ Ce^{-\gamma_1n}.
\end{equation}

\vspace{15pt}
\textbf{(B2)} There do not exist an ergodic acip $\mu_k$ of $F$, $\nu\in\bR\setminus\{0\}$, $c\in S^1$ and a measurable $\phi:\bT^2\to S^1$ such that
\[
e^{i\nu\tau}=c\,\frac{\phi\circ F}{\phi}\qquad\mu_k\text{-a.e.}
\]

\begin{oss}
Condition \textbf{(B2)} is the (multiplicative) non-arithmeticity condition. It contains both the non-coboundary condition (if $\tau=\Psi\circ F-\Psi+c$ with $\Psi$ measurable, take $\phi=e^{i\nu\Psi}$) and the non-lattice condition (if $\tau+\Psi-\Psi\circ F\in a\bZ+b$, take $\nu=2\pi/a$). It is essential that $\phi$ be allowed to take values in $S^1$ (see Remark \ref{rmk:counterex}).
\end{oss}

\begin{oss}[Meaning of \textbf{(B1)}]\label{rmk:B1star}
{For $v\in\fC_u\setminus\{0\}$, set
$u=D_{\frh(q)}F^nv/\|D_{\frh(q)}F^nv\|$. Then
\[
\frac{D_{(\frh(q),z)}\tilde F^n(v,0)}{\|D_{\frh(q)}F^nv\|}
=\big(u,d(\tau_n\circ\frh)_q(u)\big).
\]
Thus \textbf{(B1)} bounds the total inverse-Jacobian weight of branches whose lifted expanding cone contains a direction close to the plane $\{(v,t):\zeta(v)=t\}$. This is a plane-transversality condition for the three-dimensional extension.

The cotangent map of $\tilde F$ sends $(\eta,\nu)$ to
$((D_pF^*)^{-1}(\eta-\nu d\tau_p),\nu)$. On covectors $(\zeta,-1)$, which annihilate the above planes, it induces
\[
G(p,\zeta)=\big(F(p),(D_pF^*)^{-1}(\zeta+d\tau_p)\big).
\]
If $G^n(p,\zeta_0)=(q,\zeta)$ and $p=\frh(q)$, then
$\zeta-d(\tau_n\circ\frh)_q=(D_q\frh)^*\zeta_0$.
Consequently, branches with $\|\zeta_0\|\le R$ are counted in \textbf{(B1)} whenever $C_\star R\le e^{\kappa n}$, by Lemma~\ref{lem:apriori}(b). This implication is the captivity mechanism used in Section~\ref{sec:proof-main2} and analogous to the partial captivity condition in \cite{F11,NTW}.
}\end{oss}

\begin{oss}[On the hypotheses on $\tau$]\label{rmk:counterex}
Let $F(x,\theta)=(\ell x,\theta+\omega(x))$ with $\ell\in\bN$ large and $\omega$ smooth and not cohomologous to a constant, and let $\tau=\omega$. Then $\tau$ satisfies the pairwise condition
\[
|\partial_x(\tau_n\circ\frh)-\partial_x(\tau_n\circ\frh')|\ge \chi_u\,(J_\frh+J_{\frh'})
\]
for all transversal pairs of inverse branches, as well as the non-coboundary and non-lattice conditions with real continuous $\Psi$. However, since $\omega=\theta\circ F-\theta\pmod 1$, the function $e^{2\pi i(z-\theta)}$ is $\tilde F$-invariant and $\cL_{2\pi k}=e^{2\pi ik\theta}\cL_Fe^{-2\pi ik\theta}$ has the eigenvalue $1$ for each $k\in\bZ$; hence the conclusions of Theorems \ref{thm:DI} and \ref{teoB} fail. Accordingly, \textbf{(B2)} fails (with $\phi=e^{2\pi ik\theta}$) and so does \textbf{(B1)} (see Lemma \ref{lem:skew} with $b=1$). {In particular, non-coboundary and non-lattice conditions formulated only with continuous real-valued transfer functions,} or pairwise lower bounds on the derivatives of $\tau_n$ along transversal inverse branches, do not suffice. In Appendix \ref{app:B1} we characterise \textbf{(B1)} in the case of skew products.

\end{oss}

\begin{oss}[Sharpness]\label{rmk:skew-sharp}
In the setting of Lemma \ref{lem:skew}, condition (iii) is also necessary for Theorem \ref{thm:DI}. Indeed, suppose $\tau-b\omega=C_0\circ f-C_0+k$. Functions of the form $e^{2\pi im\theta}u(x)$ are mapped by $\cL_\nu$ to $e^{2\pi im\theta}\cL_f\big(e^{i\nu\tau-2\pi im\omega}u\big)$. If $b\neq0$, take $\nu_m:=2\pi m/b$; then $e^{i\nu_m\tau-2\pi im\omega}=e^{i\nu_m(C_0\circ f-C_0+k)}$, {and, taking the smooth representative of $C_0$, the function $u=e^{i\nu_mC_0}h$ belongs to every $H^s$ and gives}
\[
\cL_{\nu_m}\big(e^{2\pi im\theta}e^{i\nu_mC_0}h\big)=e^{i\nu_mk}\,e^{2\pi im\theta}e^{i\nu_mC_0}h .
\]
Hence $\cL_{\nu_m}$ has an eigenvalue of modulus one for a sequence $|\nu_m|\to\infty$, which is incompatible with \eqref{main-est}. If $b=0$ the same holds with $m=0$ and every $\nu$. {Under the VE hypotheses of Theorem~\ref{thm:DI}, this proves the equivalence of \textbf{(B1)} and its Dolgopyat conclusion for skew products with $\tau=\tau(x)$.}
\end{oss}

{
\subsection{Limit theorems and higher order corrections}\label{sec:limit-definitions}
We study the statistical properties of the Birkhoff sums
\[
\tau_n=\sum_{j=0}^{n-1}\tau\circ F^j.
\]
}

For this application we assume that $\mu_k(\tau)=0$ for every ergodic acip $\mu_k$, $1\le k\le\ell$. This is a componentwise assumption: subtracting a single constant achieves it only when all component means agree. The multiplicative non-arithmeticity assumption is \textbf{(B2)} and it excludes observables cohomologous to lattice-valued functions. We do not consider the lattice case here and in particular, the equality $\cL_{\nu+2\pi k}=\cL_\nu$ holds for integer-valued observables.

Write
\[
\fn_\sigma(x)=\frac{1}{\sqrt{2\pi\sigma^2}}e^{-x^2/(2\sigma^2)},
\qquad \Phi_\sigma(z)=\int_{-\infty}^z\fn_\sigma(x)\,\dd x.
\]
Let $\Pi u=\sum_{k=1}^{\ell}\rho_k\ell_k(u)$ be the spectral projector of $\cL_F$ at $1$, with $\int\rho_k\,\dd\frm=1$. Here $\ell_k$ are positive invariant linear functionals, normalized by $\ell_k(\rho_j)=\delta_{kj}$ and $\sum_k\ell_k(u)=\int u\,\dd\frm$. For an initial probability density $v\in H^s$, put $c_k(v)=\ell_k(v)$. Then $c_k(v)\ge0$ and $\sum_kc_k(v)=1$.

The Edgeworth expansion of order $R$ means
\begin{equation}\label{GenEdgeExp}
\bP_{v\frm}\!\left(\frac{\tau_n}{\sqrt n}\le z\right)
=\sum_{k=1}^{\ell}c_k(v)\Phi_{\sigma_k}(z)
+\sum_{k=1}^{\ell}\sum_{j=1}^{R}n^{-j/2}P_{j,k,v}(z)\fn_{\sigma_k}(z)
+o(n^{-R/2}),
\end{equation}
uniformly in $z\in\bR$. The variances $\sigma_k^2>0$ are defined in Lemma~\ref{lem:branches} below. For an ergodic invariant initial measure only its component occurs in the leading term. In general the limit is a mixture of Gaussians. The order-one expansion gives a Berry--Esseen bound relative to this mixture, with a constant depending on $F$, $\tau$, and the initial density.

For endpoint functions $\psi\in H^s$ and $\xi\in L^2$, we use the notation
\[
I_n(g;\psi,\xi):=\int_{\bT^2}\psi\,g(\tau_n)\,\xi\circ F^n\,\dd\frm,
\qquad c_k(\psi,\xi):=\ell_k(\psi)\int\xi\rho_k\,\dd\frm.
\]
Here $\psi$ incorporates both the initial density and the
initial observable. More precisely, for an initial probability
measure $v\,\frm$, with $v\in H^s$, and an observable
$a\in\cC^{\lceil s\rceil}$, we set $\psi=av$.
The Sobolev multiplication estimate gives
\[
\|\psi\|_{H^s}
\le C\|a\|_{\cC^{\lceil s\rceil}}\|v\|_{H^s},
\]
so the theorem applies. More generally, it suffices that
multiplication by $a$ be bounded on $H^s$, and we do not require a Banach algebra property.

The mixing global expansion of order $R$ will have the form
\begin{equation}\label{eq:global-expansion}
I_n(g;\psi,\xi)
=\sum_{j=0}^{R}n^{-(j+1)/2}
\int_{\bR}K_{j,\psi,\xi}(x/\sqrt n)g(x)\,\dd x
+o(n^{-(R+1)/2}),
\end{equation}
where $K_{j,\psi,\xi}=\sum_k R_{j,k,\psi,\xi}\fn_{\sigma_k}$ and the $R_{j,k,\psi,\xi}$ are polynomials of degree at most $3j$. In particular,
\[
K_{0,\psi,\xi}(x)=\sum_kc_k(\psi,\xi)\fn_{\sigma_k}(x).
\]
The mixing local expansion of order $R$ means
\begin{equation}\label{eq:local-expansion}
\sqrt n\,I_n(g;\psi,\xi)
=\sum_{m=0}^{\lfloor R/2\rfloor}n^{-m}
\int_{\bR}Q_{m,\psi,\xi}(x)g(x)\,\dd x+o(n^{-R/2}),
\end{equation}
where $\deg Q_{m,\psi,\xi}\le2m$ and
\[
Q_{0,\psi,\xi}=\sum_k\frac{c_k(\psi,\xi)}{\sqrt{2\pi}\sigma_k}.
\]
The coefficients depend on the initial and terminal functions through the spectral amplitudes, as well as on the cumulants of $\tau_n$. In the stationary ergodic case, $\psi=a\rho_k$ gives $c_k(\psi,\xi)=\mu_k(a)\mu_k(\xi)$ whenever $a\rho_k\in H^s$.

{
\subsection{Function spaces for the mixing expansions}
We use the Fourier transform
\begin{equation}\label{hat g}
\widehat g(\nu)=\int_{\bR}e^{-i\nu x}g(x)\,\dd x.
\end{equation}
}

For $p>1$ and an integer $k\ge0$, set
\[
C^p(g)=\sup_{\nu\in\bR}(1+|\nu|)^p|\widehat g(\nu)|,
\quad M_k(g)=\int_{\bR}(1+|x|^k)|g(x)|\,\dd x,
\quad \|g\|_{\fF_k^p}=C^p(g)+M_k(g).
\]
We write $g\in\fF_k^p$ when $g$ is continuous and this norm is finite. In particular, $\widehat g\in\cC^k$ and
\[
\max_{0\le j\le k}\|\widehat g^{(j)}\|_\infty\le M_k(g).
\]
The explicit absolute-moment requirement ensures that the polynomial integrals in \eqref{eq:local-expansion} are ordinary Lebesgue integrals. Compactly supported smooth functions belong to all these spaces. For integer $p$, integrable derivatives through order $p$ imply $C^p(g)<\infty$ by integration by parts. We use $p=q+2$ with any fixed $q>0$ in Theorem~\ref{teoA}.

\subsection{Main Results}\label{sec:MainResults}

For $\nu\in\bR$ let $\cL_\nu u:=\cL_F(e^{i\nu\tau}u)$ be the twisted transfer operator, where $\cL_F$ is the transfer operator of $F$ (see Section \ref{sec:transfer} for the precise definitions), and for $s>0$ let
\[
\|u\|_\nu:=\|u\|_{L^2}+(1+|\nu|)^{-s}\|u\|_{H^s}.
\]
Let $\sigma_F\in(0,1)$ be the rate of the Lasota--Yorke inequality for $\cL_F$ on $H^s$, given by Corollary \ref{thm:LY-F} below, and let $\Theta$, $\gamma_1$, $\kappa$ be as in \eqref{eq:Linfty-growth} and \textbf{(B1)}. We will assume the following compatibility condition:

\vspace{10pt}
\noindent\textbf{(P)}\quad $\Theta\,\sigma_F^2<1,\qquad \log\Theta<\gamma_1,\qquad 2\log\Theta<\kappa s.$
\vspace{10pt}

\begin{oss}
Condition \textbf{(P)} requires the possible growth of
$\|\cL_F^n1\|_\infty$ to be sufficiently slow compared with
the Sobolev contraction and the cancellation supplied by
\textbf{(B1)}. More precisely, it reads
\[
\log\Theta<
\min\left\{
2\log(\sigma_F^{-1}),\,
\gamma_1,\,
\frac{\kappa s}{2}
\right\}.
\]
The first inequality allows the Sobolev contraction to absorb
the weak-norm growth, while the remaining inequalities ensure
that the cancellation and frequency truncation estimates
still yield contraction after that growth is taken into account.

This restriction is natural for the systems considered here,
whose central direction is neutral or close to neutral.
For skew products one can take $\Theta=1$, and
\textbf{(P)} is automatic. More generally, the bound
$\Theta=\lambda_c$ shows that it suffices for the central
rate $\lambda_c$ to be sufficiently close to $1$ relative
to the other rates. In a family of systems, this conclusion
requires uniform bounds keeping $\sigma_F$ below $1$ and
$\gamma_1,\kappa s$ above $0$.
{Thus \textbf{(P)} compares central growth with the available contraction and transversality rates; it does not follow from partial hyperbolicity alone. Moreover, if $s>1$, Sobolev embedding and \eqref{LY-F} give
\[
\|\cL_F^n1\|_\infty\le C\|\cL_F^n1\|_{H^s}\le C.
\]
In that case one may take $\Theta=1$, and \textbf{(P)} is again automatic.}
\end{oss}

{Our principal technical result is the following Dolgopyat-type estimate, proved in Section~\ref{sec:proof-main}.}

\begin{thm}[Dolgopyat inequality]\label{thm:DI}
Let $s>0$ be such that $F$ is $s$--VE and satisfies \eqref{eq:Linfty-growth}, let $\tau \in \cC^{\infty}(\bT^2, \bR)$ satisfy \textbf{(B1)}, and assume \textbf{(P)}. There exist $\nu_0>1$, $C>0$, $\gamma>0$ and $\delta\in(0,1)$ such that\footnote{With a slight abuse of notation we write $\|\cL\|_{\nu}$ for the operator norm of $\cL$ acting on $H^{s}$ endowed with $\|\cdot\|_\nu$.}
\begin{equation}\label{main-est}
\|\cL^{n}_\nu \|_{\nu}\le C|\nu|^{\delta} e^{-\gamma n}, \qquad \text{for all} \quad n\in\bN,\ |\nu|\ge \nu_0.
\end{equation}
In particular, $\|\cL^{n}_\nu \|_{H^s\to H^s}\le C|\nu|^{s+\delta} e^{-\gamma n}$ for all $n\in\bN$ and $|\nu|\ge\nu_0$.
\end{thm}

The above Theorem~\ref{thm:DI} is the key to prove the following Theorems~\ref{teoA} and~\ref{teoB}.

\begin{bigthm}[Limit theorems and higher order expansions]\label{teoA}
Let $F$ be an $s$--VE SVPH satisfying \textbf{(A6)}, and let $\tau\in\cC^\infty(\bT^2,\bR)$ satisfy \textbf{(B1)}, \textbf{(B2)}, and $\mu_k(\tau)=0$ for every ergodic acip $\mu_k$. Assume \textbf{(P)}. Then all component variances $\sigma_k^2$ are strictly positive. For every integer $R\ge1$:
\begin{enumerate}[label=(\alph*),leftmargin=*]
\item Every initial probability measure $v\frm$ with $v\in H^s$ admits the Edgeworth expansion \eqref{GenEdgeExp}, uniformly in $z$.
\item For every $q>0$, $g\in\fF_0^{q+2}$, $\psi\in H^s$, and $\xi\in H^s$, the mixing global expansion \eqref{eq:global-expansion} holds.
\item For every $q>0$, $g\in\fF_{R+1}^{q+2}$, $\psi\in H^s$, and $\xi\in H^s$, the mixing local expansion \eqref{eq:local-expansion} holds.
\end{enumerate}
In (b)--(c) one may more generally take $\xi\in L^2$. The coefficients are independent of $g$; the remainders are understood for fixed $v$, or fixed $g,\psi,\xi$, respectively.
\end{bigthm}
The proof of the above theorem is given in Section \ref{sec:proof-thmA}.
\begin{oss}
The spectral gap and perturbation near zero frequency provide the cumulants and the Gaussian terms. The additional role of the Dolgopyat estimate is to control large real frequencies, allowing expansions of arbitrary order. The general spectral methods underlying these expansions are developed in \cite{FeLi,FePe}; Section~\ref{sec:proof-thmA} gives the required argument for the present, possibly non-ergodic, Sobolev setting.
\end{oss}

The second application concerns the exponential decay of correlations for the skew product $\tilde F\in \cC^\infty(\bT^3,\bT^3)$ given by
\begin{equation}\label{F*}
\tilde F(x,y,z)=(F(x,y), z+\tau(x,y)) \quad \mod \bZ^3.
\end{equation}
It is proved in Section \ref{sec:doc}.

\begin{bigthm}[Exponential mixing]\label{teoB}
Let $s>0$ be such that $F$ is $s$--VE and satisfies \textbf{(A6)}, let $\tau \in \cC^{\infty}(\bT^2, \bR)$ satisfy \textbf{(B1)} and \textbf{(B2)}, assume \textbf{(P)}, and let ${\mu}_*$ be an ergodic acip for $F$ with density $h_*\in H^s$.\footnote{$H^s$ denotes the standard Sobolev space, with norm \eqref{equivalence of norm H L}. See Section \ref{sec:transfer} for definitions and results on the invariant measures for the map $F$.} Then the $\tilde F$-invariant absolutely continuous measure ${\dd\mu(p,z)=h_*(p)\,\dd p\,\dd z}$ is exponentially mixing with respect to H\"older observables; that is, for each $\alpha>0$ there exist $\gamma>0$ and $C>0$ such that, for any $f,g\in \cC^\alpha(\bT^3,\bC)$,
\[
\left|\int_{\bT^3} f\cdot g\circ \tilde F^n \, \dd{\mu}-\int_{\bT^3} f \, \dd{\mu} \int_{\bT^3} g \, \dd{\mu} \right|\le C \|f\|_{\cC^\alpha}\|g\|_{\cC^\alpha} e^{-\gamma n}, \quad \forall n \ge 1.
\]
\end{bigthm}

\begin{oss}\label{rmk-teoB}{
The circle extension is related at the operator level to a suspension semiflow over $F$. To define a roof one must require positivity; for a general real $\tau$, one may instead use $r=\tau+c>0$. On the imaginary axis this constant shift only multiplies the twisted operators by scalars of modulus one. For a positive roof $\tau$ (see, e.g., \cite{PoYu, Tsu08}),} the Fourier modes in $z$ of the correlations of $\tilde F$ are governed by the operators $\cL_{2\pi k}$, which are the restriction to the imaginary axis of the operators $\cL_F(e^{-z\tau}\,\cdot\,)$, $z\in\bC$, governing the correlations of the semiflow. As recalled in the introduction, in the skew-product case $\tilde F(x,\theta,z)=(f(x),\theta+\omega(x),z+\tau(x))$ exponential mixing already follows from \cite{ChHu}. Theorem \ref{teoB} applies also when the base $F$ is a partially expanding map which is not a skew product, or when $\tau$ depends on $\theta$, provided the hypotheses are satisfied.
\end{oss}

\section{Transfer operators}\label{sec:transfer}
{This section collects the transfer-operator tools used in the Dolgopyat estimate and both applications. The twisted Lasota--Yorke inequality is proved in Section~\ref{subsec:LYnu}.}

\subsection{Spectral gap}

The Ruelle transfer operator associated to $F$ is defined via
\begin{equation}\label{L_F}
    \cL_F u(x)=\sum_{y\in F^{-1}(x)} \frac{u(y)}{|\det DF(y)|}, \quad u\in L^1(\bT^2), x,y \in \bT^2,
\end{equation}
and used to study the statistical properties of $F$. The main result of \cite{Tsu1} is the quasi-compactness of $\cL_F$ on suitable Sobolev spaces, which will be crucial for our arguments. To state the result, let us recall the definition of the spaces used: for each $s\in \bR^+$, $H^s(\bT^2, \bC)$ will denote the standard Hilbert space of $L^2(\bT^2, \bC)$ functions endowed with the norm\footnote{Here $\hat{\cdot}$ denotes the Fourier transform and $\langle k \rangle= \sqrt{1+\|k\|^2}$.}
\begin{equation}
\label{equivalence of norm H L}
\|u\|_{H^s}:= \left( \sum_{k \in \bZ^2} \langle k \rangle^{2s} |\hat{u}(k)|^2  \right)^{\frac 12}.
\end{equation}

We first relate \eqref{virtually} to the condition used in \cite{Tsu1}.

\begin{lem}\label{lem:VE-Tsujii}
Let $F$ be $s$-VE. Then, for every $t\in\big[\frac s2,s\big]$,
\[
\lim_{n\to\infty}B_{2t}(F^n)^{1/n}<1.
\]
\end{lem}

\begin{proof}
If $(F^n)^\dagger(p,\eta)=(q,\xi)$, then $\xi=[(D_pF^n)^*]^{-1}\eta$, hence, by \eqref{DF-1}, $\|\xi\|/\|\eta\|\le\|(D_pF^n)^{-1}\|\le C\lambda_c^n$. Since $0\le 2t-s\le s$ and $\lambda_c\ge1$,
\[
\left(\frac{\|\xi\|}{\|\eta\|}\right)^{2t}\le (C\lambda_c^n)^{2t-s}\left(\frac{\|\xi\|}{\|\eta\|}\right)^{s}\le C\lambda_c^{sn}\left(\frac{\|\xi\|}{\|\eta\|}\right)^{s}.
\]
Thus $B_{2t}(F^n)\le C\lambda_c^{sn}B_s(F^n)$, and the claim follows from \eqref{virtually}.
\end{proof}

\begin{thm}[\cite{Tsu1}]\label{thm:Tsujii}
Let $t>0$ be such that $\lim_{n\to\infty}B_{2t}(F^n)^{1/n}<1$. Then $\cL_F$ extends to a bounded operator on $H^t(\bT^2)$, which is quasi-compact with essential spectral radius strictly smaller than $1$.
\end{thm}

\begin{cor}[DFLY inequality for \texorpdfstring{$\cL_F$}{Lg}]\label{thm:LY-F}
Let $F$ be $s$-VE for some $s>0$. For any $s'\in \bR$ such that $\max\big\{\frac s2, s- \frac{1}{2}\big\}\le s' <s$ there exist $\sigma_* \in (0,1)$ and $C>0$ such that:
\begin{equation}\label{LY-F}
\begin{split}
&\|\cL_F^n u\|_{H^{s'}}\le C\|u\|_{H^{s'}}\\
&\|\cL_F^n u\|_{H^s}\le C\sigma_*^n \|u\|_{H^s}+C \|u\|_{H^{s'}}
\end{split}
\end{equation}
for each $n \in \bN$ and $u\in H^s(\bT^2,\bC)$.
\end{cor}

\begin{proof}
Let $t\in\{s',s\}$. By Lemma \ref{lem:VE-Tsujii} and Theorem \ref{thm:Tsujii}, $\cL_F$ is quasi-compact on $H^t$ with essential spectral radius $\vartheta_t<1$. Since $H^t\subset L^1$ and $\|\cL_F^nu\|_{L^1}\le\|u\|_{L^1}$, every eigenvalue $z$ of $\cL_F$ on $H^t$ satisfies $|z|\le1$. Moreover, eigenvalues of modulus one have no Jordan blocks: if $\cL_Fu=zu$ and $\cL_Fv=zv+u$ with $|z|=1$, then $\cL_F^nv=z^nv+nz^{n-1}u$, whence $n\|u\|_{L^1}\le2\|v\|_{L^1}$ for all $n$, i.e.\ $u=0$. It follows that $\sup_n\|\cL_F^n\|_{H^t}<\infty$, which for $t=s'$ is the first inequality.

Choose $\sigma_*\in(\max\{\vartheta_s,\vartheta_{s'}\},1)$ such that $\cL_F$ has no spectrum on $\{|z|=\sigma_*\}$, neither on $H^s$ nor on $H^{s'}$, and let $P$, $P'$ be the spectral projectors of $\cL_F$ on $H^s$, $H^{s'}$, respectively, associated with the spectrum in $\{|z|>\sigma_*\}$; they have finite rank. For $|z|$ large the resolvents of $\cL_F$ on $H^s$ and on $H^{s'}$ are given by the same Neumann series, so they coincide on $H^s$; by analytic continuation they coincide on $H^s$ for all $z$ in $\{|z|\ge\sigma_*\}$ outside a finite set, hence $P'|_{H^s}=P$. Since $H^s$ is dense in $H^{s'}$ and $P'$ is continuous with finite rank, $P'(H^{s'})=P(H^s)\subset H^s$. On this finite dimensional space all norms are equivalent, so 
$$\|Pu\|_{H^s}=\|P'u\|_{H^s}\le C\|P'u\|_{H^{s'}}\le C\|u\|_{H^{s'}}$$
 for $u\in H^s$. Finally, 
 $$\|\cL_F^n(\Id-P)\|_{H^s}\le C\sigma_*^n,$$ and therefore
\[
\|\cL_F^nu\|_{H^s}\le\|\cL_F^n(\Id-P)u\|_{H^s}+\|\cL_F^n\|_{H^s}\|Pu\|_{H^s}\le C\sigma_*^n\|u\|_{H^s}+C\|u\|_{H^{s'}}. \qedhere
\]
\end{proof}

{Fix a number $\sigma_F\in(0,1)$ for which \eqref{LY-F} holds for some admissible $s'$. It may be chosen arbitrarily close to the infimum of the admissible rates; we do not assume that this infimum is positive or attained.}

Next, we define twisted transfer operators.

\begin{defn}
For any $\nu \in \bR$ and $u\in L^1(\bT^2, \bC)$ we define the twisted transfer operator $\cL_\nu$ associated to $F$ by
\[
\cL_\nu u (x,y)= \sum_{(w_1,w_2)\in F^{-1}(x,y)} \dfrac{u(w_1,w_2)}{\det DF(w_1,w_2)} e^{i\nu \tau(w_1,w_2)}.
\]
\end{defn}

Iterating the above formula we see that for each $n\in \bN$
\begin{equation}\label{L^n}
\cL^n_\nu u (x,y)=\sum_{(w_1,w_2)\in F^{-n}(x,y)} \dfrac{u(w_1,w_2)}{\det DF^n(w_1,w_2)} e^{i\nu \tau_n(w_1,w_2)},
\end{equation}
where $\tau_n:=\sum_{k=0}^{n-1}\tau \circ F^k$. Note also that
\begin{equation}\label{eq:TTO}
\cL_\nu f=\cL_F(e^{i\nu\tau } f).
\end{equation}

Finally, let us define the norm we will be using. For each $\nu\in \bR$ and $s\in \bR^+$ let
\[
\|u\|_{\nu}:=\|u\|_{L^2}+(1+|\nu|)^{-s}\|u\|_{H^s}.
\]
For each fixed $\nu$, the above norm is equivalent to $\|\cdot \|_{H^s}$ (not uniformly in $\nu$).

The following is a Lasota-Yorke type inequality with strong norm $\|\cdot \|_{\nu}$ and weak norm $\|\cdot\|_{L^2}$, whose proof is postponed to Section \ref{subsec:LYnu}.

\begin{lem}[DFLY inequality for \texorpdfstring{$\cL_\nu$}{Lg}]\label{lem:LYnu}
Let $F$ be $s$--VE for some $s>0$, satisfying \eqref{eq:Linfty-growth}, and $\tau\in\cC^\infty$. For each $\sigma \in (\sigma_F,1)$ there exist $C_\sigma>0$ and $C_{\tau}>0$ such that, for all $n\in\bN$, $u\in H^s$, $\nu \in \bR$,
\begin{equation}\label{LYnu}
\|\cL_\nu^n u\|_{\nu}\le C_{\sigma}\sigma^n\|u\|_{\nu}+C_{\tau}\Theta^{n/2}\| u \|_{L^{2}}.
\end{equation}
\end{lem}

\begin{oss}
{The proof of Theorem~\ref{thm:DI} specifies admissible $\gamma$ and $\delta$ in terms of the Lasota--Yorke constants, the rates in \eqref{eq:Linfty-growth} and \textbf{(B1)}, and the inverse-branch derivative bounds of Lemma~\ref{lem:apriori}. These bounds also depend on finitely many derivatives of $F$ and $\tau$.} 
\end{oss}

\subsection{Ergodic properties}\label{sec:ErgodicOp}

Let $\frm$ denote normalised Lebesgue measure on $\bT^2$.
We record the ergodic consequences of quasi-compactness
that will be used below. In particular, we distinguish the
invariant probabilities from the dual eigenfunctionals
of the transfer operator.

\begin{prop}\label{prop:ergodic-decomposition}
Suppose that $F$ is $s$--VE for some $s>0$. Then:
\begin{enumerate}[label=(\roman*),leftmargin=*]
\item
The operator $\cL_F:H^s\to H^s$ has spectral radius one
and essential spectral radius strictly smaller than one.
Its peripheral eigenvalues are semisimple, and
\[
\sup_{n\ge0}\|\cL_F^n\|_{H^s\to H^s}<\infty.
\]
Moreover, $\cL_F$ is positive and
$\|\cL_F^n\|_{L^1\to L^1}\le1$.

\item
There are nonnegative functions
$\rho_1,\dots,\rho_\ell\in H^s$, satisfying
\[
\int\rho_k\,\dd\frm=1,
\qquad
\rho_j\rho_k=0\quad\frm\text{-a.e.}\quad(j\ne k),
\]
which form a basis of $\ker(\cL_F-\Id)$ on $H^s$.
The probabilities $\mu_k=\rho_k\frm$ are precisely the
ergodic absolutely continuous invariant probabilities.
Every absolutely continuous invariant probability is a
convex combination of the $\mu_k$.

\item
Let $\Pi$ denote the spectral projector associated with
the eigenvalue $1$. There are positive bounded
functionals $\ell_k\in(H^s)'$ such that
\[
\Pi u=\sum_{k=1}^{\ell}\rho_k\ell_k(u),
\]
and
\[
\ell_k\cL_F=\ell_k,\qquad
\ell_k(\rho_j)=\delta_{kj},\qquad
\sum_{k=1}^{\ell}\ell_k(u)=\int u\,\dd\frm.
\]
Here positivity means $\ell_k(u)\ge0$ whenever
$u\in H^s$ is real and nonnegative almost everywhere.
Furthermore,
\[
\frac1N\sum_{n=0}^{N-1}\cL_F^n
\longrightarrow\Pi
\quad\text{in operator norm on }H^s.
\]

\item
The physical basins of $\mu_1,\dots,\mu_\ell$ cover
$\frm$-almost every point.\end{enumerate}
\end{prop}

\begin{proof}
Write $L=\cL_F$.
Positivity, preservation of integrals and
$L^1$-contractivity follow from the inverse-branch formula:
\[
|Lu|\le L|u|,
\qquad
\int Lu\,\dd\frm=\int u\,\dd\frm.
\]
By Lemma~\ref{lem:VE-Tsujii} and
Theorem~\ref{thm:Tsujii}, the essential spectral radius
of $L$ on $H^s$ is strictly smaller than one.
Since $H^s$ embeds continuously into $L^1$, contractivity
excludes eigenvalues of modulus greater than one.

It also excludes nontrivial Jordan blocks at peripheral
eigenvalues. Indeed, if
\[
Lu=zu,\qquad Lv=zv+u,\qquad |z|=1,
\]
then
\[
L^nv=z^nv+nz^{n-1}u.
\]
Consequently,
\[
n\|u\|_{L^1}
\le \|L^nv\|_{L^1}+\|v\|_{L^1}
\le2\|v\|_{L^1},
\]
so $u=0$.
Finally, integration is a nonzero invariant bounded
functional on $H^s$. Hence $1$ belongs to the spectrum
of $L$, and, since it lies outside the essential
spectrum, it is an eigenvalue.
These facts imply that the spectral radius is one and
that the powers of $L$ are uniformly bounded on $H^s$.

Let
\[
A_N=\frac1N\sum_{n=0}^{N-1}L^n.
\]
The spectral decomposition, semisimplicity of the
peripheral eigenvalues, and the fact that
$N^{-1}\sum_{n=0}^{N-1}z^n\to0$ for $|z|=1$, $z\ne1$,
give
\begin{equation}\label{eq:Cesaro-projector}
A_N\longrightarrow\Pi
\quad\text{in operator norm on }H^s.
\end{equation}
In particular, $\Pi$ is positive and preserves integrals.\\
Next, set
\[
E=\ker(L-\Id)\subset H^s.
\]
This space is finite-dimensional and therefore closed
also as a subspace of $L^1$.

Let $v\in L^1$ satisfy $Lv=v$. Choose $u_j\in H^s$
with $u_j\to v$ in $L^1$. Since $A_Nv=v$ and
$\|A_N\|_{L^1\to L^1}\le1$, letting $N\to\infty$ gives
\[
\|v-\Pi u_j\|_{L^1}\le\|v-u_j\|_{L^1}.
\]
Thus $v$ belongs to the $L^1$-closure of $E$, hence to
$E$ itself. We have proved
\begin{equation}\label{eq:fixed-L1-Hs}
\ker(L-\Id)\big|_{L^1}
=
\ker(L-\Id)\big|_{H^s}.
\end{equation}

The real fixed space $E_{\mathbb R}$ is a vector
sublattice of $L^1$. Indeed, if $u\in E_{\mathbb R}$,
then
\[
|u|=|Lu|\le L|u|.
\]
Both sides have the same integral, so $L|u|=|u|$.
By \eqref{eq:fixed-L1-Hs}, $|u|\in E_{\mathbb R}$.
The usual formulas for the maximum and minimum of two
real functions now give the sublattice property.

A finite-dimensional vector sublattice of $L^1$ has
a basis of nonnegative, pairwise disjoint functions:
its positive cone is generated by finitely many
disjoint atoms. Normalising these atoms by their
integrals gives $\rho_1,\dots,\rho_\ell$.
Every nonnegative fixed function consequently has a
unique representation
\begin{equation}\label{eq:decomp}
v=\sum_{k=1}^{\ell}c_k\rho_k,\qquad c_k\ge0.
\end{equation}
For a probability density, $\sum_kc_k=1$.

Each $\mu_k=\rho_k\frm$ is ergodic. Otherwise, an
invariant measurable set of $\mu_k$-measure strictly
between zero and one would split $\mu_k$ into two
nonzero invariant absolutely continuous measures.
Their densities would belong to $E$ by
\eqref{eq:fixed-L1-Hs}, contradicting that $\rho_k$
generates an extreme ray of its positive cone.
Conversely, \eqref{eq:decomp} shows
that every ergodic acip is one of the $\mu_k$.\\
Next, define $\ell_k$ by the unique expansion
\[
\Pi u=\sum_k\rho_k\ell_k(u).
\]
They are bounded because $\Pi$ is bounded and its range
is finite-dimensional. Their positivity follows from
positivity of $\Pi$ and disjointness of the $\rho_k$.
The identities
\[
L\Pi=\Pi L=\Pi,\qquad \Pi\rho_j=\rho_j
\]
give $\ell_kL=\ell_k$ and
$\ell_k(\rho_j)=\delta_{kj}$.
Since $\Pi$ preserves integrals,
\[
\sum_k\ell_k(u)=\int\Pi u\,\dd\frm=\int u\,\dd\frm.
\]

Finally, for each $k$, let
\[
\mathcal B_k=
\left\{
p:\frac1N\sum_{n=0}^{N-1}\varphi(F^np)
\longrightarrow\int\varphi\,\dd\mu_k
\text{ for every }\varphi\in C^0(\bT^2)
\right\}.
\]
Using a countable uniformly dense subset of
$C^0(\bT^2)$, these are measurable sets.
Birkhoff's theorem and ergodicity give
$\mu_k(\mathcal B_k)=1$.

 Note that 
$F^{-1}\mathcal B_k=\mathcal B_k$.
Set
\[
B=\bT^2\setminus\bigcup_k\mathcal B_k.
\]
Then $F^{-1}B=B$ and $\mu_k(B)=0$ for every $k$.
Transfer-operator duality gives
\[
\int_B A_N1\,\dd\frm
=\frac1N\sum_{n=0}^{N-1}
  \int 1_B\circ F^n\,\dd\frm
=\frm(B).
\]
By \eqref{eq:Cesaro-projector}, $A_N1\to\Pi1$ in
$L^1$. Since $\Pi1$ is a linear combination of the
$\rho_k$, it follows that
\[
\frm(B)=\int_B\Pi1\,\dd\frm=0.
\]
Moreover, $\mu_k(\mathcal B_k)=1$ and $\mu_k\ll\frm$
imply $\frm(\mathcal B_k)>0$, so each $\mu_k$ is
physical.
\end{proof}

\begin{oss}\label{rmk:dual-functionals}
The functionals $\ell_k$ are dual eigenfunctionals of
$\cL_F$; they should not be identified with the
invariant probabilities $\mu_k$.
In fact, for $u\in H^s$,
\[
\ell_k(u)=\int_{\mathcal B_k}u\,\dd\frm.
\]
Indeed, the basins are pairwise disjoint,
$F^{-1}\mathcal B_k=\mathcal B_k$, and
$\mu_j(\mathcal B_k)=\delta_{jk}$. Integrating
$A_Nu\to\Pi u$ over $\mathcal B_k$ proves the identity.
\end{oss}

Under \textbf{(A6)}, Lemma~\ref{lem:exact} below excludes
peripheral eigenvalues other than $1$. Consequently,
there exist $C>0$ and $\vartheta\in(0,1)$ such that
\[
\|\cL_F^n-\Pi\|_{H^s\to H^s}
\le C\vartheta^n.
\]
Without \textbf{(A6)}, the Ces\`aro convergence in
Proposition~\ref{prop:ergodic-decomposition} remains
valid, but convergence of the iterates themselves is
not asserted.
Virtually expanding maps are \textit{volume expanding} (see \cite{Tsu1}), in the following sense
\begin{equation}
    \lim_{n \to \infty} \left(\inf_{p \in \bT^2}|\det D_pF^n|\right)^{1/n}>1\,.
\end{equation}

In the present notation, \textbf{(A6)} says that $(F,\mu_k)$ is weakly mixing for each $k\in\{1,\dots,\ell\}$, and \eqref{eq:Linfty-growth} reads $\|\cL_F^n1\|_{L^\infty}\le C_\infty\Theta^n$.

\begin{lem}\label{lem:exact}
Assume \textbf{(A6)}. Then there are no peripheral eigenvalues\footnote{Eigenvalues on the unit circle $\{z \in \mathbb C\mid |z|=1\}$.} of $\cL_F: H^s \to H^s$ other than $1$. 
\end{lem}

\begin{proof}
Suppose that $\cL_F\rho =e^{i\lambda}\rho$ with $\rho \in H^s\setminus\{0\}$ and $\lambda\in\bR$. Since $|\rho|=|\cL_F\rho|\le \cL_F|\rho|$ pointwise and $\int\cL_F|\rho|\,\dd\frm=\int|\rho|\,\dd\frm$, we have $\cL_F|\rho|=|\rho|$ $\frm$-a.e. Hence $|\rho|\,\frm$ is a finite $F$-invariant absolutely continuous measure, so by (ii) $|\rho|=\sum_k c_k\rho_k$ with $c_k\ge 0$, and some $c_k>0$.

Write $\rho=g|\rho|$ with $|g|=1$ on $\{|\rho|>0\}$. Then, for $\frm$-a.e.\ $x$,
\[
\Big|\sum_{y\in F^{-1}(x)}\frac{g(y)|\rho|(y)}{|\det D_yF|}\Big|=\sum_{y\in F^{-1}(x)}\frac{|\rho|(y)}{|\det D_yF|},
\]
and equality in the triangle inequality, together with $\cL_F\rho=e^{i\lambda}g\,\cL_F|\rho|$, yields $g(y)=e^{i\lambda}g(F(y))$ for all $y\in F^{-1}(x)$ with $|\rho|(y)>0$. Thus $g\circ F=e^{-i\lambda}g$ $\mu_k$-a.e., i.e.\ $g\in L^2(\mu_k)$ is a non-zero eigenfunction of the Koopman operator of $(F,\mu_k)$ with eigenvalue $e^{-i\lambda}$. Since $(F,\mu_k)$ is weakly mixing, $e^{-i\lambda}=1$.
\end{proof}

\subsection{Spectral properties of the twisted operators}

\begin{lem}\label{lem:spectral}
Let $F$ be $s$-VE and satisfy \textbf{(A6)}, and let $\tau \in \cC^{\infty}$ satisfy \textbf{(B2)}. Then:
\begin{enumerate}[label=(\alph*)]
\item $1$ is the only peripheral eigenvalue of $\cL_F$ on $H^s$; the corresponding spectral projector $\Pi$ satisfies $\cL_F\Pi=\Pi\cL_F=\Pi$ and $\|\cL_F^n(\Id-\Pi)\|_{H^s}\le C\theta_0^n$ for some $C>0$, $\theta_0\in(0,1)$;
\item for all $\nu \in \bR\setminus\{0\}$, the spectral radius of $\cL_\nu$ on $H^s$ is smaller than $1$;
\item for each compact $J\subset\bR\setminus\{0\}$ there exist $C_J>0$ and $\theta_J\in(0,1)$ such that
\begin{equation}\label{midOpEst}
\sup_{\nu\in J}\|\cL^n_\nu\|_{H^s} \le C_J\,\theta_J^n\qquad\forall n\in\bN .
\end{equation}
\end{enumerate}
\end{lem}

\begin{proof}
(a) By Proposition \ref{prop:ergodic-decomposition}(i) and Lemma \ref{lem:exact}, $1$ is the only peripheral eigenvalue and the rest of the spectrum is contained in a disc of radius $\theta_0<1$. Since $\sup_n\|\cL_F^n\|_{H^s}<\infty$ by \eqref{LY-F}, the eigenvalue $1$ is semisimple, i.e.\ $\cL_F\Pi=\Pi$.

(b) {Fix $\sigma\in(\sigma_F,1)$. Lemma~\ref{lem:LYnu} and the compact embedding $H^s\hookrightarrow L^2$ imply that the essential spectral radius of $\cL_\nu$ on $H^s$ is at most $\sigma<1$ \cite{HeHe}. Consequently any spectral value on or outside the unit circle would be an eigenvalue.} Since $\|\cL_\nu^nu\|_{L^1}\le\|\cL^n_F|u|\|_{L^1}=\|u\|_{L^1}$, every eigenvalue $z$ satisfies $|z|\le1$. Suppose $\cL_\nu u=e^{i\vartheta}u$ with $u\in H^s\setminus\{0\}$. Then $|u|\le\cL_F|u|$ and $\int\cL_F|u|=\int|u|$, so $\cL_F|u|=|u|$ a.e.; by Proposition \ref{prop:ergodic-decomposition}(ii), $|u|=\sum_kc_k\rho_k$ with some $c_k>0$. Writing $u=g|u|$ with $|g|=1$ on $\{|u|>0\}$, the identity
\[
\sum_{y\in F^{-1}(x)}\frac{e^{i\nu\tau(y)}g(y)|u|(y)}{|\det D_yF|}=e^{i\vartheta}g(x)\sum_{y\in F^{-1}(x)}\frac{|u|(y)}{|\det D_yF|}
\]
forces equality in the triangle inequality, i.e.\ $e^{i\nu\tau(y)}g(y)=e^{i\vartheta}g(F(y))$ for $\mu_k$-a.e.\ $y$. Thus $e^{i\nu\tau}=e^{i\vartheta}\,g\circ F/g$ $\mu_k$-a.e., contradicting \textbf{(B2)}.

(c) Since $\tau\in\cC^\infty$, multiplication by $\tau$ is bounded on $H^s$, so $\cL_\nu=\sum_{k\ge0}\frac{\nu^k}{k!}\cL_F((i\tau)^k\,\cdot\,)$ converges in $L(H^s,H^s)$ and $\nu\mapsto\cL_\nu$ is continuous in operator norm. By (b), for each $\nu\in J$ there is $m_\nu\in\bN$ with $\|\cL_\nu^{m_\nu}\|_{H^s}\le\frac12$; by continuity the same holds, with $\frac12$ replaced by $\frac34$, for $\nu'$ in a neighbourhood of $\nu$. Covering $J$ by finitely many such neighbourhoods, and using $\sup_{\nu\in J,\,m\le M}\|\cL^m_\nu\|_{H^s}<\infty$ for each $M$, we obtain \eqref{midOpEst}.
\end{proof}

\section{Proof of the Dolgopyat estimate (Theorem \ref{thm:DI})}\label{sec:proof-main}

This section is devoted to the proof of Theorem \ref{thm:DI}. Estimates of this type go back to \cite{Dima98, Dima02}; for twisted transfer operators of partially expanding maps see also \cite{F11, BuEs}. Throughout, $C$ denotes a generic positive constant (with possible subscripts indicating dependence) that does not depend on $n\in\bN$ or $\nu\in \bR$, and may change from line to line. We fix $\rho>0$, to be chosen in Section \ref{sec:proof-main2}, and set
\[
n(\nu):=\lceil\rho\log|\nu|\rceil,\qquad |\nu|>1.
\]

\subsection*{Road map to the proof} The proof follows the logic of \cite{BuEs} and it has two main steps. First, we reduce the Dolgopyat
estimate to an $L^2$ contraction estimate at a time proportional
to $\log|\nu|$. We then prove this contraction by localising
in space and frequency and applying the transversality
assumption \eqref{eq:B1star}.

Section~\ref{subsec:LYnu} establishes the twisted
Lasota--Yorke inequality of Lemma~\ref{lem:LYnu}, using
the multiplication estimate in Lemma~\ref{lem:twist-mult}.
This separates two cases. When the Sobolev part of
$\|h\|_\nu$ is large relative to its $L^2$ part,
Lemma~\ref{lem:case1} gives contraction after a fixed
number of iterates. For the remaining functions,
Corollary~\ref{key-cor} shows that it suffices to obtain
the $L^2$ estimate \eqref{L^2est} at the logarithmic time
$n(\nu)$. The same corollary explains how these two
contraction mechanisms combine to yield the estimate
for arbitrary times.

Section~\ref{sec:proof-main2} proves the required
$L^2$ estimate. After establishing bounds on inverse
branches and a non-stationary phase estimate, we truncate
the input frequencies and localise in space. For each
inverse branch, we distinguish the output frequencies
compatible with the phase gradient from those outside
this region. Lemma~\ref{lem:captivity} uses
\eqref{eq:B1star} to bound the total Jacobian weight of
branches compatible with any given frequency;
Lemma~\ref{lem:trapped} converts this bound into an
$L^2$ estimate. The complementary contribution is
controlled by non-stationary phase in
Lemma~\ref{lem:untrapped}. Finally,
Proposition~\ref{prop:L2est} combines these estimates
with the frequency truncation error and verifies the
hypotheses of Corollary~\ref{key-cor}.

\subsection{Proof of Lemma \ref{lem:LYnu}}\label{subsec:LYnu}
We begin with a classic result.
\begin{lem}[Sobolev multiplication]\label{lem:leib}
For $s>0$ there is $C_s>0$ such that
\[
\|fu\|_{H^s}\le C_s\|f\|_{\cC^{\lceil s\rceil}}\|u\|_{H^s}
\]
for $f\in\cC^{\lceil s\rceil}$ and $u\in H^s$.
\end{lem}
\begin{proof}
For integer $s$ this follows from the Leibniz rule. For noninteger $s$, apply that rule through order $\lfloor s\rfloor$ and use the fractional difference norm at the remaining order; see also \cite{Babook}.
\end{proof}
\begin{lem}[Twisted multiplication]\label{lem:twist-mult}
There exists $C_0>0$ and, for each $N\in\bN$, a constant $C_N>0$ such that, for all $\nu\in\bR$ and $u\in H^s$,
\[
\|e^{i\nu\tau_N}u\|_{H^s}\le C_0\|u\|_{H^s}+C_N(1+|\nu|)^s\|u\|_{L^2}.
\]
\end{lem}

\begin{proof}
Let $f:=e^{i\nu\tau_N}$. Then $\|fu\|_{L^2}=\|u\|_{L^2}$ and, for each $j\ge1$, $\|D^jf\|_{\cC^0}\le C_N(1+|\nu|)^j$. We will repeatedly use the interpolation inequality: for $0\le t<s$ and $\ve>0$,
\begin{equation}\label{eq:interp-nu}
(1+|\nu|)^{s-t}\|u\|_{H^t}\le(1+|\nu|)^{s-t}\|u\|_{L^2}^{1-\frac ts}\|u\|_{H^s}^{\frac ts}\le\ve\|u\|_{H^s}+C_\ve(1+|\nu|)^s\|u\|_{L^2},
\end{equation}
which follows from Young's inequality. Write $s=k+\sigma$ with $k\in\bN\cup\{0\}$ and $\sigma\in(0,1]$. We use the equivalent norm $\|u\|_{L^2}+\sum_{|\beta|=k}\|\partial^\beta u\|_{H^\sigma}$, where, for $\sigma\in(0,1)$, $\|w\|_{H^\sigma}$ is equivalent to $\|w\|_{L^2}+[w]_\sigma$, with $[w]_\sigma^2:=\iint\frac{|w(x)-w(y)|^2}{d(x,y)^{2+2\sigma}}\dd x\,\dd y$.

Let $g\in\cC^1$, $w\in H^\sigma$. If $\sigma\in(0,1)$, from $$|g(x)w(x)-g(y)w(y)|\le|g(x)||w(x)-w(y)|+|g(x)-g(y)||w(y)|,$$ we have 
$$|g(x)-g(y)|\le\min\{2\|g\|_{\cC^0},\|Dg\|_{\cC^0}d(x,y)\}$$ and from
$$\int\min\{a^2,b^2d(x,y)^2\}d(x,y)^{-2-2\sigma}\dd x\le Ca^{2-2\sigma}b^{2\sigma}$$, we get
\begin{equation}\label{eq:prod-sigma}
[gw]_\sigma\le\|g\|_{\cC^0}[w]_\sigma+C\|g\|_{\cC^0}^{1-\sigma}\|Dg\|_{\cC^0}^{\sigma}\|w\|_{L^2}.
\end{equation}
For $\sigma=1$ the same bound follows from $\nabla(gw)=g\nabla w+w\nabla g$.\\
For $|\beta|=k$, by the Leibniz rule $\partial^\beta(fu)=f\partial^\beta u+\sum_{0<\gamma\le\beta}c_{\beta,\gamma}\,\partial^\gamma f\,\partial^{\beta-\gamma}u$. By \eqref{eq:prod-sigma} with $g=f$, $|f|=1$,
\[
\|f\partial^\beta u\|_{H^\sigma}\le C\|\partial^\beta u\|_{H^\sigma}+C_N(1+|\nu|)^\sigma\|u\|_{H^k},
\]
with $C$ depending only on $s$.

For the remaining terms we can argue as follows. For $0<\gamma\le\beta$, $|\gamma|=j\ge1$, by \eqref{eq:prod-sigma} with $g=\partial^\gamma f$,
\[
\|\partial^\gamma f\,\partial^{\beta-\gamma}u\|_{H^\sigma}\le C_N\big((1+|\nu|)^j\|u\|_{H^{s-j}}+(1+|\nu|)^{j+\sigma}\|u\|_{H^{k-j}}\big).
\]
Since $(s-j)+j=(k-j)+(j+\sigma)=k+\sigma=s$ and $k+\sigma=s$, all the terms with a factor $C_N$ are of the form $(1+|\nu|)^{s-t}\|u\|_{H^t}$ with $t<s$, and by \eqref{eq:interp-nu} with $\ve=(C_N)^{-1}$ they are bounded by $\|u\|_{H^s}+C_N'(1+|\nu|)^s\|u\|_{L^2}$. Summing over $\beta$ gives the claim, with $C_0$ depending only on $s$.
\end{proof}

Note that, by Cauchy--Schwarz, $\big|\cL_\nu^{n}u\big|^2\le \cL_F^n1\cdot\cL_F^n|u|^2$; hence, by \eqref{eq:Linfty-growth},
\begin{equation}\label{eq:L2-bound}
\|\cL^n_\nu u\|_{L^2}\le \|\cL^n_F1\|_{L^\infty}^{1/2}\|u\|_{L^2}\le C_\infty^{1/2}\Theta^{n/2}\|u\|_{L^2}\qquad\forall n\in\bN,\ \nu\in\bR.
\end{equation}

\begin{proof}[Proof of Lemma \ref{lem:LYnu}]
{Let $\sigma\in(\sigma_F,1)$ and choose $\sigma_*\in(\sigma_F,\sigma)$. Inequality \eqref{LY-F} holds with this larger rate and the admissible $s'$ chosen above.} Fix $N\in\bN$. By \eqref{eq:TTO}, \eqref{LY-F} and Lemma \ref{lem:twist-mult} (applied with $s$ and with $s'$),
\[
\|\cL_\nu^N u\|_{H^s} \le C\sigma_*^N\big(C_0\|u\|_{H^s}+C_N(1+|\nu|)^s\|u\|_{L^2}\big)+C\big(C_0\|u\|_{H^{s'}}+C_N(1+|\nu|)^{s'}\|u\|_{L^2}\big).
\]
By interpolation, for each $\ve>0$, $\|u\|_{H^{s'}}\le\ve\|u\|_{H^s}+C_\ve\|u\|_{L^2}$. Multiplying by $(1+|\nu|)^{-s}\le1$,
\[
(1+|\nu|)^{-s}\|\cL_\nu^N u\|_{H^s}\le (CC_0\sigma_*^N+CC_0\ve)(1+|\nu|)^{-s}\|u\|_{H^s}+C_{N,\ve}\|u\|_{L^2}.
\]
Since $\sigma_*<\sigma$, we can choose $N$ so large that $CC_0\sigma_*^N\le\frac12\sigma^N$, and then $\ve$ so small that $CC_0\ve\le\frac12\sigma^N$. Together with \eqref{eq:L2-bound},
\[
\|\cL_\nu^N u\|_{\nu}\le \sigma^N\|u\|_\nu+C'\|u\|_{L^2}.
\]
Applying this to $\cL^{kN}_\nu u$ and using \eqref{eq:L2-bound}, we get $X_{k+1}\le\sigma^NX_k+C''\Theta^{kN/2}\|u\|_{L^2}$ for $X_k:=\|\cL_\nu^{kN} u\|_{\nu}$. Since $\sigma<1\le\Theta^{1/2}$, induction gives $X_k\le\sigma^{kN}\|u\|_\nu+C'''\Theta^{kN/2}\|u\|_{L^2}$. The same computation with $N$ replaced by $m<N$ gives $\|\cL_\nu^m u\|_\nu\le C\|u\|_\nu$. Writing $n=kN+m$ yields \eqref{LYnu}.
\end{proof}

\begin{lem}\label{lem:case1}
Let $\sigma$, $C_\sigma$ and $C_{\tau}$ be as in Lemma \ref{lem:LYnu}. There exist $n_0:=n_0(\sigma)\in \bN$ and $\sigma_1\in (0,1)$ such that, setting $\overline C:= \max\{C_\sigma, C_{\tau}\Theta^{n_0/2}\}$, for each $\nu\in \bR$, if $h\in H^{s}$ satisfies
\[
2\overline C (1+|\nu|)^{s} \|h\|_{L^2}\le \|h\|_{H^{s}},
\]
then $\|\cL_\nu^{n_0}h\|_{\nu} \le \sigma_1\|h\|_{\nu}$.
\end{lem}

\begin{proof}
Fix $\sigma_1\in(\tfrac12,1)$ and choose $n_0\in\bN$ such that $C_\sigma \sigma^{n_0} \leq \sigma_1 - \frac 12$. By Lemma \ref{lem:LYnu},
\[
\|\cL_\nu^{n_0}h\|_{\nu}\le C_\sigma \sigma^{n_0}\|h\|_{\nu}+C_{\tau}\Theta^{n_0/2}\|h\|_{L^2}.
\]
Using the assumption on $h$ and $\overline C\ge C_{\tau}\Theta^{n_0/2}$,
\[
C_{\tau}\Theta^{n_0/2}\|h\|_{L^2}\le \overline C\,\|h\|_{L^2}\le \tfrac12 (1+|\nu|)^{-s}\|h\|_{H^{s}}\le \tfrac12 \|h\|_{\nu}.
\]
Therefore,
\[
\|\cL_\nu^{n_0}h\|_{\nu}\le \left(C_\sigma\sigma^{n_0}+\tfrac12\right)\|h\|_{\nu}\le \sigma_1\|h\|_{\nu}. \qedhere
\]
\end{proof}

\begin{cor}[$L^2$-estimate reduction]\label{key-cor}
Fix $\sigma\in(\sigma_F,1)$ with $\sigma\Theta^{1/2}<1$ in Lemma \ref{lem:LYnu} (possible by \textbf{(P)}), and assume that $\rho\log\Theta<\frac14$. Assume there exist $\nu_0>1$, $C_1>0$ and $\gamma_4>\frac12\log\Theta$ such that, for all $|\nu|\ge\nu_0$ and all $h \in H^s$ with $\|h\|_{H^{s}}\le 2 \overline C (1+|\nu|)^{s} \|h\|_{L^2}$,
\begin{equation}\label{L^2est}
\|\cL^{n(\nu)}_\nu h \|_{L^2}\le C_1 e^{-n(\nu)\gamma_4}\|h\|_{L^2}.
\end{equation}
Then, up to enlarging $\nu_0$, Theorem \ref{thm:DI} holds, i.e.\ there exist $C>0$, $\gamma>0$ and $\delta\in(0,1)$ such that \eqref{main-est} holds.
\end{cor}

\begin{proof}
Fix $|\nu|\ge \nu_0$ and set $N:=n(\nu)$; enlarging $\nu_0$ we may assume $2N\ge n_0$. By Lemma \ref{lem:LYnu} and $\|\cdot\|_{L^2}\le\|\cdot\|_\nu$, $\|\cL^N_\nu h\|_\nu\le C\Theta^{N/2}\|h\|_\nu$; applying Lemma \ref{lem:LYnu} once more,
\begin{equation}\label{eq:key-cor}
\|\cL^{2N}_\nu h\|_{\nu}\le C\sigma^N \|\cL^N_\nu h\|_{\nu}+C\Theta^{N/2}\|\cL_\nu^N h\|_{L^2}\le C(\sigma\Theta^{1/2})^N \|h\|_{\nu}+C\Theta^{N/2}\|\cL_\nu^N h\|_{L^2}.
\end{equation}
\emph{Case 1:} $\|h\|_{H^s}\ge 2\overline C(1+|\nu|)^s\|h\|_{L^2}$. By Lemma \ref{lem:case1}, $\|\cL_\nu^{n_0}h\|_\nu\le e^{-\gamma_6n_0}\|h\|_\nu$ with $\gamma_6:=\frac{\log\sigma_1^{-1}}{n_0}$.

\noindent\emph{Case 2:} $\|h\|_{H^s}\le 2\overline C(1+|\nu|)^s\|h\|_{L^2}$. By \eqref{eq:key-cor} and \eqref{L^2est}, $\|\cL^{2N}_\nu h\|_{\nu}\le C\big((\sigma\Theta^{1/2})^N+C_1e^{-(\gamma_4-\frac12\log\Theta)N}\big)\|h\|_\nu\le e^{-2\gamma_5N}\|h\|_\nu$, with $\gamma_5:=\frac14\min\{\log(\sigma\Theta^{1/2})^{-1},\gamma_4-\frac12\log\Theta\}>0$, up to enlarging $\nu_0$.

Let $\gamma'=\min\{\gamma_5,\gamma_6\}$ and $\gamma:=\min\{\gamma',(4\rho)^{-1}\}$. We have shown that for each $h\in H^s$ there is $t(h)\in\{n_0,2N\}$ with $\|\cL^{t(h)}_\nu h\|_\nu\le e^{-\gamma' t(h)}\|h\|_\nu$. Set $h_0=h$, $h_{j+1}=\cL_\nu^{t(h_j)}h_j$ and $T_j=\sum_{i<j}t(h_i)$, so that $\|h_j\|_\nu\le e^{-\gamma' T_j}\|h\|_\nu\le e^{-\gamma T_j}\|h\|_\nu$. Given $n$, let $j$ be such that $T_j\le n<T_{j+1}$; then $n-T_j<2N$, so $e^{-\gamma T_j}\le e^{2\gamma N}e^{-\gamma n}$ and, since $\sup_{m<2N}\|\cL^m_\nu\|_\nu\le C\Theta^{N}\le C|\nu|^{\rho\log\Theta}$ by \eqref{LYnu},
\[
\|\cL^n_\nu h\|_\nu=\|\cL^{n-T_j}_\nu h_j\|_\nu\le C|\nu|^{\rho\log\Theta}e^{2\gamma N}e^{-\gamma n}\|h\|_\nu\le C|\nu|^{\rho\log\Theta+2\gamma\rho}e^{-\gamma n}\|h\|_\nu.
\]
The claim follows with $\delta:=\rho\log\Theta+2\gamma\rho<\frac34$.
\end{proof}

Thanks to Corollary \ref{key-cor}, we are left to prove \eqref{L^2est} for functions satisfying $\|h\|_{H^{s}}\le 2\overline C (1+|\nu|)^{s} \|h\|_{L^2}$. This is done in the next subsection.

\subsection{Proof of \eqref{L^2est}: a captivity estimate}\label{sec:proof-main2}

\subsubsection{Preliminaries}

For $q\in\bT^2$ and $n\in\bN$ we denote by $\cH^n(q)$ the set of inverse branches of $F^n$ defined on a ball $B$ centred at $q$ of radius smaller than $\frac14$ (they exist since $F^n:\bT^2\to\bT^2$ is a covering map and $B$ is simply connected). For $\frh\in\cH^n(q)$ we write
\[
J_\frh(y):=|\det D_{\frh(y)}F^n|^{-1},\qquad c_\frh(y):=d\big(\tau_n\circ\frh\big)_y\in(\bR^2)^*,\qquad y\in B .
\]
With this notation, for $y\in B$,
\begin{equation}\label{eq:Ln-branches}
\cL_\nu^n u(y)=\sum_{\frh\in\cH^n(q)}J_\frh(y)\,e^{i\nu\tau_n(\frh(y))}\,u(\frh(y)).
\end{equation}

\begin{lem}[A priori bounds]\label{lem:apriori}
For $j\ge1$ and a map $\varphi$ defined on an open set of $\bT^2$ (with values in $\bT^2$, lifted locally to $\bR^2$) let $|\varphi|_{j}:=\max_{1\le i\le j}\sup\|D^i\varphi\|$; for a real function $\psi$ let $\|\psi\|_{\cC^j}:=\sup|\psi|+|\psi|_j$. There exist constants $C_1,C_2,\dots>0$ and $p_1\le p_2\le\dots$, with $p_1\ge\log(\Lambda\lambda_c\Theta)$, such that, for all $n\in\bN$, $q\in\bT^2$, $\frh\in\cH^n(q)$ (defined on the ball $B$ centred at $q$):
\begin{enumerate}[label=(\alph*)]
\item $\sup_B\|D\frh\|\le C_1\lambda_c^n$;
\item if $u\in D_{\frh(y)}F^n\fC_u$, $y\in B$, is a unit vector, then $\|D_y\frh\,u\|\le C_\star\lambda^{-n}$;
\item $|\frh|_j+\|\tau_n\circ\frh\|_{\cC^j}+\|\log J_\frh\|_{\cC^j}\le C_je^{p_jn}$ and $\|D^jJ_\frh\|\le C_jJ_\frh\,e^{jp_jn}$ on $B$;
\item $\sharp\cH^n(q)\le C_\infty C_*(\Lambda\lambda_c\Theta)^n$.
\end{enumerate}
\end{lem}

\begin{proof} We divide the proof in steps.\\
\emph{Preliminaries.} Since $F^n$ is a local diffeomorphism of a compact manifold, it is a covering map; as $B$ is simply connected, $F^{-n}(B)$ is a disjoint union of open sets, each mapped diffeomorphically onto $B$ by $F^n$, and the elements of $\cH^n(q)$ are the inverses of these restrictions. In particular $\frh\mapsto\frh(q)$ is a bijection between $\cH^n(q)$ and $F^{-n}(q)$. For $\frh\in\cH^n(q)$ and $0\le k\le n$ set
\[
\frh^{(k)}:=F^{n-k}\circ\frh:B\to\bT^2 .
\]
Since $F^k\circ\frh^{(k)}=F^n\circ\frh=\Id_B$, we have $\frh^{(k)}\in\cH^k(q)$, $\frh^{(0)}=\Id_B$ and $\frh^{(n)}=\frh$. Moreover $F$ is injective on $\frh^{(k+1)}(B)$ (because $F\circ\frh^{(k+1)}=\frh^{(k)}$ is injective), so $g_k:=\big(F|_{\frh^{(k+1)}(B)}\big)^{-1}:\frh^{(k)}(B)\to\frh^{(k+1)}(B)$ is well defined and
\begin{equation}\label{eq:apriori-rec}
\frh^{(k+1)}=g_k\circ\frh^{(k)},\qquad Dg_k=(DF)^{-1}\circ g_k .
\end{equation}
By \textbf{(A1)} and compactness, $\inf_{\bT^2}\det DF>0$, so $\sup\|(DF)^{-1}\|<\infty$; differentiating $Dg_k=(DF)^{-1}\circ g_k$ repeatedly, $D^ig_k$ is a polynomial in $(DF)^{-1}\circ g_k$ and $D^lF\circ g_k$, $l\le i$. Hence there are constants $G_i$, depending only on $\|F\|_{\cC^{i}}$ and $\sup\|(DF)^{-1}\|$, with
\begin{equation}\label{eq:apriori-G}
\sup\|D^ig_k\|\le G_i\qquad\text{for all }i\ge1,\ k\ge0 .
\end{equation}

\emph{(a).} Let $p\in\bT^2$, $k\in\bN$ and $w\in\cT_{F^k(p)}\bT^2$. Since $\chi_u\chi_c<1$, there is $\theta_0>0$, depending only on $\chi_u,\chi_c$, such that the angle between any nonzero vector of $\fC_u$ and any nonzero vector of $\fC_c$ is at least $\theta_0$. By \eqref{invariance of cone}, $D_pF^k\fC_u\subset\fC_u$; hence we can write $w=e+c$ with $e\in D_pF^k\fC_u$, $c\in\fC_c$ and $\|e\|,\|c\|\le\|w\|/\sin\theta_0$. Write $e=D_pF^kv$ with $v\in\fC_u$; since $\fC_u\setminus\{0\}\subset\bR^2\setminus\fC_c$, \eqref{partial hyperbolicity 2} gives $\|e\|\ge C_\star^{-1}\lambda^k\|v\|$, i.e.\ $\|(D_pF^k)^{-1}e\|\le C_\star\lambda^{-k}\|e\|$. Moreover $\|(D_pF^k)^{-1}c\|\le\lambda^+_{c,k}(p)\|c\|\le C_\star(\lambda_c^+)^k\|c\|$. Since $\lambda>1$ and $\lambda_c^+\le\lambda_c$, $\lambda_c\ge1$,
\begin{equation}\label{eq:apriori-a}
\|(D_pF^k)^{-1}\|\le\frac{2C_\star}{\sin\theta_0}\,\lambda_c^k=:C_1'\lambda_c^k .
\end{equation}
As $D_y\frh^{(k)}=(D_{\frh^{(k)}(y)}F^k)^{-1}$, this gives $\sup_B\|D\frh^{(k)}\|\le C_1'\lambda_c^k$ for all $0\le k\le n$, and (a) for $k=n$.

\emph{(b).} Write $u=D_{\frh(y)}F^nv/\|D_{\frh(y)}F^nv\|$ with $v\in\fC_u$. As above, $\|D_{\frh(y)}F^nv\|\ge C_\star^{-1}\lambda^n\|v\|$, hence $\|D_y\frh\,u\|=\|v\|/\|D_{\frh(y)}F^nv\|\le C_\star\lambda^{-n}$.

\emph{(c), bounds on $\frh$.} Let $a_k^{(i)}:=\sup_B\|D^i\frh^{(k)}\|$. We prove by induction on $i$ that there exist $c_i>0$ and $\pi_i$ with
\begin{equation}\label{eq:apriori-ind}
a_k^{(i)}\le c_ie^{\pi_ik}\qquad\text{for all }k\ge0,
\end{equation}
{Enlarge the bounds $G_i$ so that $G_i\ge1$, and choose $c_i$ nondecreasing.} Here $\pi_1:=\log(\Lambda\lambda_c\Theta)$ and $\pi_i:=i\,\pi_{i-1}+\log G_1+1$ for $i\ge2$. For $i=1$ this is \eqref{eq:apriori-a}, since $\lambda_c\le\Lambda\lambda_c$. Let $i\ge2$ and assume \eqref{eq:apriori-ind} for $1,\dots,i-1$. By \eqref{eq:apriori-rec} and the Faà di Bruno formula,
\[
D^i\frh^{(k+1)}=Dg_k(\frh^{(k)})\,D^i\frh^{(k)}+\sum_{l=2}^{i}\ \sum_{\substack{b_1+\dots+b_l=i\\ b_1,\dots,b_l\ge1}}c_{i,b}\,D^lg_k(\frh^{(k)})\big[D^{b_1}\frh^{(k)},\dots,D^{b_l}\frh^{(k)}\big],
\]
with combinatorial coefficients $c_{i,b}$. In the second sum $l\ge2$, so each $b_m\le i-1$, and by the induction hypothesis and $\pi_1\le\dots\le\pi_{i-1}$ each product is bounded by $c_{i-1}^ie^{l\pi_{i-1}k}\le c_{i-1}^ie^{i\pi_{i-1}k}$. Using \eqref{eq:apriori-G},
\[
a_{k+1}^{(i)}\le G_1a_k^{(i)}+K_ie^{i\pi_{i-1}k},\qquad a_0^{(i)}=0,
\]
with $K_i$ depending only on $i$, $G_1,\dots,G_i$, $c_{i-1}$. Iterating,
\[
a_k^{(i)}\le K_i\sum_{m=0}^{k-1}G_1^{k-1-m}e^{i\pi_{i-1}m}\le K_i\,k\,e^{(i\pi_{i-1}+\log G_1)k}\le c_ie^{\pi_ik},
\]
which proves \eqref{eq:apriori-ind}. Note that $\pi_i$ is nondecreasing in $i$ (as $\pi_{i-1}>0$) and $|\frh|_j=\max_{i\le j}a^{(i)}_n\le c_je^{\pi_jn}$.

\emph{(c), bounds on $\tau_n\circ\frh$ and $\log J_\frh$.} Since $F^m\circ\frh=\frh^{(n-m)}$,
\[
\tau_n\circ\frh=\sum_{m=0}^{n-1}\tau\circ F^m\circ\frh=\sum_{k=1}^{n}\tau\circ\frh^{(k)},\qquad
\log J_\frh=-\sum_{m=0}^{n-1}\log\det DF\circ F^m\circ\frh=\sum_{k=1}^{n}\psi\circ\frh^{(k)},
\]
where $\psi:=-\log\det DF\in\cC^\infty(\bT^2)$ by \textbf{(A1)}. Let $\phi\in\{\tau,\psi\}$. Then $\sup|\phi\circ\frh^{(k)}|\le\|\phi\|_{\cC^0}$ and, by Faà di Bruno and \eqref{eq:apriori-ind}, for $1\le i\le j$,
\[
\sup_B\|D^i(\phi\circ\frh^{(k)})\|\le\sum_{l=1}^{i}\sum_{b_1+\dots+b_l=i}c_{i,b}\|D^l\phi\|_{\cC^0}\prod_{m=1}^la_k^{(b_m)}\le K'_j\|\phi\|_{\cC^j}e^{j\pi_jk}.
\]
Summing over $1\le k\le n$ gives $\|\tau_n\circ\frh\|_{\cC^j}+\|\log J_\frh\|_{\cC^j}\le K''_j\,n\,e^{j\pi_jn}$. Setting
\[
p_j:=j\,\pi_j+1,
\]
which is nondecreasing in $j$ and satisfies $p_1=\pi_1+1\ge\log(\Lambda\lambda_c\Theta)$, and choosing $C_j$ large, we obtain the first estimate in (c). For the second, $D^jJ_\frh=D^j\big(e^{\log J_\frh}\big)$ is $J_\frh$ times a polynomial (a complete Bell polynomial) in $D\log J_\frh,\dots,D^j\log J_\frh$, which is homogeneous of degree $j$ when $D^i$ is given weight $i$; hence $\|D^jJ_\frh\|\le C_jJ_\frh\,e^{jp_jn}$.

\emph{(d).} By \eqref{det-est}, $|\det D_{\frh(q)}F^n|\le C_*\Lambda^n\lambda_c^n$, i.e.\ $1\le C_*(\Lambda\lambda_c)^nJ_\frh(q)$. Summing over $\frh\in\cH^n(q)$ and using that $\frh\mapsto\frh(q)$ is a bijection onto $F^{-n}(q)$,
\[
\sharp\cH^n(q)\le C_*(\Lambda\lambda_c)^n\sum_{y\in F^{-n}(q)}|\det D_yF^n|^{-1}\le C_\infty C_*(\Lambda\lambda_c\Theta)^n
\]
by \eqref{eq:Linfty-growth}.
\end{proof}

\begin{lem}[Non-stationary phase]\label{lem:nsp}
Let $m\in\bN$. There exists $C_m>0$ such that the following holds. Let $U\subset\bR^2$ be open, $r\in(0,1]$, $a\in\cC^m_c(U)$ with $\|D^ja\|_{\cC^0}\le A\,r^{-j}$ for $0\le j\le m$, and $\Phi\in\cC^{m+1}(U,\bR)$ with $\|D^j\Phi\|_{\cC^0}\le Q$ for $2\le j\le m+1$, where $Q\ge1$, and $|\nabla\Phi|\ge\mu\ge1$ on $\operatorname{supp}a$. Then
\[
{\Big|\int_U e^{i\Phi}a\Big|\le C_m\,A\,\operatorname{Leb}(\operatorname{supp}a)}\,\Big(\frac{1}{r\mu}+\frac{Q}{\mu^2}\Big)^m.
\]
\end{lem}

\begin{proof}
Set $V=\nabla\Phi/(i|\nabla\Phi|^2)$ on a neighbourhood of $\supp a$, so that $V\cdot\nabla e^{i\Phi}=e^{i\Phi}$. A cutoff allows an arbitrary smooth extension away from this neighbourhood. Write $L^Tb=-\operatorname{div}(bV)$. Repeated integration by parts gives
\[
\int_Ue^{i\Phi}a=\int_Ue^{i\Phi}(L^T)^ma.
\]
On $\supp a$, the chain rule gives
\[
|V|\le\mu^{-1},\qquad |D^jV|\le C_j\mu^{-1}(1+Q/\mu)^j\quad(1\le j\le m).
\]
Indeed, derivatives of $v\mapsto v/|v|^2$ of order $l$ are bounded by $C_l|v|^{-1-l}$; the remaining factors are derivatives of $\Phi$ of order at least two. Expanding $(L^T)^ma$ produces products with $m$ factors consisting of derivatives of $V$, and total derivative order $m$, distributed between these factors and $a$. Hence
\[
|(L^T)^ma|\le C_mA\mu^{-m}\sum_{j=0}^m r^{-j}(1+Q/\mu)^{m-j}
\le C_mA\left(\frac1{r\mu}+\frac Q{\mu^2}\right)^m,
\]
where the last step uses $r\le1$. Integrating over $\supp a$ proves the claim.
\end{proof}

\subsubsection{Localisation in space and frequency}

Fix $|\nu|\ge\nu_0$, set $N:=n(\nu)$, and let $h\in H^s$ satisfy $\|h\|_{H^s}\le2\overline C(1+|\nu|)^s\|h\|_{L^2}$. We fix the scales
\[
r:=|\nu|^{-1/2},\qquad M:=|\nu|^{5/8},\qquad R:=|\nu|^{\varepsilon},
\]
with $\varepsilon\in(0,\frac1{16})$ to be chosen.

\emph{Frequency cut-off.} Write $h=h_{\le}+h_>$, where $h_\le:=\sum_{k\in\bZ^2,\ 2\pi|k|\le R|\nu|}\hat h(k)e^{2\pi ik\cdot x}$. Then $\|h_\le\|_{L^2}\le\|h\|_{L^2}$ and
\[
\|h_>\|_{L^2}\le C(R|\nu|)^{-s}\|h\|_{H^s}\le C R^{-s}\|h\|_{L^2},
\]
so that, by \eqref{eq:L2-bound},
\begin{equation}\label{eq:tail}
\|\cL^N_\nu h_>\|_{L^2}\le C\,\Theta^{N/2}|\nu|^{-\varepsilon s}\|h\|_{L^2}.
\end{equation}

\emph{Space localisation.} Let $\{\psi_\alpha\}_{\alpha\in\cA}$ be smooth functions with $\operatorname{supp}\psi_\alpha\subset\cU_\alpha:=B(z_\alpha,r)$, $\sum_\alpha\psi_\alpha^2=1$, $\|D^j\psi_\alpha\|_{\cC^0}\le C_jr^{-j}$, with intersection multiplicity of $\{\cU_\alpha\}$ bounded by $c_\sharp$ and $\sharp\cA\le Cr^{-2}=C|\nu|$. Then
\[
\|\cL^N_\nu h_\le\|^2_{L^2}=\sum_\alpha\|\psi_\alpha\cL^N_\nu h_\le\|_{L^2}^2,\qquad \psi_\alpha\cL^N_\nu h_\le=\sum_{\frh\in\cH^N(z_\alpha)}g_{\alpha,\frh},
\]
where, by \eqref{eq:Ln-branches},
\[
g_{\alpha,\frh}(y):=\psi_\alpha(y)J_\frh(y)\,e^{i\nu\tau_N(\frh(y))}\,h_\le(\frh(y)),\qquad y\in\cU_\alpha .
\]
We regard $g_{\alpha,\frh}$ as a function on $\bR^2$ (via a chart of $\cU_\alpha$) and denote by $\hat g_{\alpha,\frh}$ its Fourier transform.

\emph{Trapped frequencies.} For $\alpha\in\cA$ and $\frh\in\cH^N(z_\alpha)$ let
\[
{\begin{gathered}\Sigma_{\alpha,\frh}:=\big\{\nu\, c_\frh(y)+(D_y\frh)^*\eta:\ y\in\cU_\alpha,\ |\eta|\le R|\nu|\big\},\\ S_{\alpha,\frh}:=\{\xi\in\bR^2:\ \operatorname{dist}(\xi,\Sigma_{\alpha,\frh})<M\}.\end{gathered}}
\]
These are the frequencies that $g_{\alpha,\frh}$ can carry: the phase of $e^{i\nu\tau_N\circ\frh}e^{2\pi ik\cdot\frh}$ has gradient $\nu c_\frh+(D\frh)^*(2\pi k)$. We split $\hat g_{\alpha,\frh}=a_{\alpha,\frh}+b_{\alpha,\frh}$ with $a_{\alpha,\frh}:=\mathds 1_{S_{\alpha,\frh}}\hat g_{\alpha,\frh}$.

\subsubsection{The trapped part: captivity}

\begin{lem}[Captivity]\label{lem:captivity}
Assume \textbf{(B1)}, and that $\rho$ and $\varepsilon$ satisfy
\begin{equation}\label{eq:cond-capt}
\varepsilon<\kappa\rho,\qquad \rho\,(\log\lambda+p_2)+\varepsilon<\tfrac38 .
\end{equation}
Then, for $|\nu|$ large enough, for all $\alpha\in\cA$ and $\xi\in\bR^2$,
\[
\sum_{\frh\in\cH^N(z_\alpha):\ \xi\in S_{\alpha,\frh}}J_\frh(z_\alpha)\le Ce^{-\gamma_1N}.
\]
\end{lem}

\begin{proof}
Let $\xi\in S_{\alpha,\frh}$. Then there are $y\in\cU_\alpha$ and $\eta$ with $|\eta|\le R|\nu|$ such that $|\xi-\nu c_\frh(y)-(D_y\frh)^*\eta|<M$. Set $\zeta:=\xi/\nu$, $\zeta_0:=\eta/\nu$ (so $|\zeta_0|\le R$), and let $u\in D_{\frh(z_\alpha)}F^N\fC_u$ be a unit vector. By Lemma \ref{lem:apriori},
\[
\begin{split}
|\zeta(u)-c_\frh(z_\alpha)(u)|&\le \frac M{|\nu|}+|c_\frh(y)-c_\frh(z_\alpha)|+|\zeta_0(D_y\frh\,u)|\\
&\le |\nu|^{-3/8}+C_2e^{p_2N}r+R\big(\|D_{z_\alpha}\frh\,u\|+C_2e^{p_2N}r\big)\\
&\le |\nu|^{-3/8}+2C_2|\nu|^{\varepsilon+\rho p_2-\frac12}+C_\star R\lambda^{-N}.
\end{split}
\]
Since $N\ge\rho\log|\nu|$ and $\varepsilon<\kappa\rho$, we have $C_\star R\lambda^{-N}\le\frac12e^{\kappa N}\lambda^{-N}$ for $|\nu|$ large. Since $N\le\rho\log|\nu|+1$, $e^{\kappa N}\lambda^{-N}\ge \lambda^{-1}|\nu|^{-\rho\log\lambda}$, and by \eqref{eq:cond-capt} the first two terms are bounded by $\frac12 e^{\kappa N}\lambda^{-N}$ for $|\nu|$ large. Hence every $\frh$ in the sum belongs to the set in \eqref{eq:B1star} (with $n=N$, $q=z_\alpha$ and this $\zeta$), and the claim follows from \textbf{(B1)}.
\end{proof}

\begin{lem}\label{lem:trapped}
Under the assumptions of Lemma \ref{lem:captivity},
\[
\sum_{\alpha\in\cA}\Big\|\sum_{\frh}a_{\alpha,\frh}\Big\|^2_{L^2}\le C e^{-\gamma_1N}\|h\|^2_{L^2}.
\]
\end{lem}

\begin{proof}
By Cauchy--Schwarz, for each $\xi$,
\[
\Big|\sum_\frh a_{\alpha,\frh}(\xi)\Big|^2\le\Big(\sum_{\frh:\ \xi\in S_{\alpha,\frh}}J_\frh(z_\alpha)\Big)\Big(\sum_\frh\frac{|\hat g_{\alpha,\frh}(\xi)|^2}{J_\frh(z_\alpha)}\Big).
\]
Integrating in $\xi$, using Lemma \ref{lem:captivity} and Plancherel,
\[
\Big\|\sum_\frh a_{\alpha,\frh}\Big\|^2_{L^2}\le Ce^{-\gamma_1N}\sum_\frh\frac{\|g_{\alpha,\frh}\|^2_{L^2}}{J_\frh(z_\alpha)}.
\]
By Lemma \ref{lem:apriori}, 
$$|\log J_\frh(y)-\log J_\frh(z_\alpha)|\le C_1e^{p_1N}r\le\log2$$
 for $y\in\cU_\alpha$ and $|\nu|$ large (as $\rho p_1<\frac12$ by \eqref{eq:cond-capt}). Hence
\[
\sum_\frh\frac{\|g_{\alpha,\frh}\|^2_{L^2}}{J_\frh(z_\alpha)}\le2\int\psi_\alpha^2\sum_\frh J_\frh|h_\le\circ\frh|^2=2\int\psi_\alpha^2\,\cL^N_F(|h_\le|^2)=2\int(\psi_\alpha^2\circ F^N)|h_\le|^2.
\]
Summing over $\alpha$ and using $\sum_\alpha\psi_\alpha^2=1$ gives the claim.
\end{proof}

\subsubsection{The untrapped part: non-stationary phase}

\begin{lem}\label{lem:untrapped}
Let $m:=64$ and assume, in addition to \eqref{eq:cond-capt}, that
\begin{equation}\label{eq:cond-nsp}
\varepsilon+\rho\, p_{m+1}\le\tfrac18 .
\end{equation}
Then, for $|\nu|$ large enough,
\[
\sum_{\alpha\in\cA}\Big\|\sum_{\frh}b_{\alpha,\frh}\Big\|^2_{L^2}\le C|\nu|^{-4}\|h\|^2_{L^2}.
\]
\end{lem}

\begin{proof}
Write $h_\le=\sum_k\hat h(k)e^{2\pi ik\cdot x}$, the sum running over $k\in\bZ^2$ with $2\pi|k|\le R|\nu|$; there are at most $C(R|\nu|)^2$ such $k$. Lifting $\frh$ to $\bR^2$,
\[
{\begin{gathered}\hat g_{\alpha,\frh}(\xi)=\sum_k\hat h(k)\,I_k(\xi),\qquad I_k(\xi):=\int e^{i\Phi_k(y)}\,\psi_\alpha(y)J_\frh(y)\,\dd y,\\ \Phi_k:=\nu\,\tau_N\circ\frh+2\pi k\cdot\frh-\xi\cdot y .\end{gathered}}
\]
We have $\nabla\Phi_k(y)=\nu c_\frh(y)+(D_y\frh)^*(2\pi k)-\xi$, hence $|\nabla\Phi_k|\ge d(\xi):=\operatorname{dist}(\xi,\Sigma_{\alpha,\frh})$ on $\cU_\alpha$, and $d(\xi)\ge M$ for $\xi\notin S_{\alpha,\frh}$. By Lemma \ref{lem:apriori}, the amplitude $\psi_\alpha J_\frh$ satisfies the assumptions of Lemma \ref{lem:nsp} with {$A=C_mJ_\frh(z_\alpha)\le C_mC_*$: bounded distortion and \eqref{eq:cond-nsp} imply, for $0\le j\le m$,
\[
\|D^j(\psi_\alpha J_\frh)\|_\infty\le C_mJ_\frh(z_\alpha)r^{-j},
\]
since $e^{p_jN}\le C|\nu|^{1/8}\le Cr^{-1}$. The bound on $A$ follows from \eqref{det-est} and $\lambda_c<\lambda$.} {Moreover,} 
$$\|D^j\Phi_k\|\le Q:=C(1+R)|\nu|e^{p_{m+1}N}$$ for $2\le j\le m+1$. For $t\ge M$,
\[
\frac1{rt}+\frac Q{t^2}\le\frac Mt\Big(\frac1{rM}+\frac Q{M^2}\Big)\le\frac Mt\,\big(|\nu|^{-1/8}+C|\nu|^{\varepsilon+\rho p_{m+1}-\frac14}\big)\le C\frac Mt|\nu|^{-1/8},
\]
by \eqref{eq:cond-nsp}. Lemma \ref{lem:nsp} gives $$|I_k(\xi)|\le C r^2\big(M|\nu|^{-1/8}/d(\xi)\big)^m$$ for $\xi\notin S_{\alpha,\frh}$. The set $\Sigma_{\alpha,\frh}$ has diameter $$D_\Sigma\le C R|\nu|e^{p_1N}\le C|\nu|^2,$$ so $$\operatorname{Leb}(\{d<t\})\le C(D_\Sigma+t)^2.$$ Decomposing $\{d\ge M\}$ into the shells $\{2^jM\le d<2^{j+1}M\}$, $j\ge0$, on which $|I_k|^2\le Cr^4|\nu|^{-m/4}2^{-2mj}$,
\[
{\begin{aligned}\int_{\{d\ge M\}}|I_k|^2\dd\xi&\le Cr^4|\nu|^{-m/4}\sum_{j\ge0}2^{-2mj}\big(D_\Sigma+2^{j+1}M\big)^2\\ &\le Cr^4|\nu|^{-m/4}\big(D_\Sigma^2+M^2\big)\le C|\nu|^{-m/4+4}.\end{aligned}}
\]
By Cauchy--Schwarz, $\sum_k|\hat h(k)|\le CR|\nu|\,\|h\|_{L^2}$, therefore $$\|b_{\alpha,\frh}\|_{L^2}\le\sum_k|\hat h(k)|\,\|I_k\mathds 1_{\{d\ge M\}}\|_{L^2}\le C R|\nu|\,|\nu|^{2-m/8}\|h\|_{L^2}.$$ Finally, by Lemma \ref{lem:apriori}(d), since $\rho p_1\le1$, we have $$\sharp\cH^N(z_\alpha)\le C(\Lambda\lambda_c\Theta)^N\le C|\nu|^{\rho\log(\Lambda\lambda_c\Theta)}\le C|\nu|$$ so
\[
\sum_\alpha\Big\|\sum_\frh b_{\alpha,\frh}\Big\|^2_{L^2}\le\sharp\cA\cdot\big(\sharp\cH^N\big)^2\max_{\alpha,\frh}\|b_{\alpha,\frh}\|^2_{L^2}\le C|\nu|^{3}\,|\nu|^{2+2\varepsilon+4-m/4}\|h\|_{L^2}^2\le C|\nu|^{-4}\|h\|^2_{L^2}. \qedhere
\]
\end{proof}

\subsubsection{Conclusion}

\begin{prop}\label{prop:L2est}
Assume \eqref{eq:Linfty-growth}, \textbf{(B1)} and \textbf{(P)}. Choose $\varepsilon:=\kappa\rho/2$ and $\rho>0$ so small that \eqref{eq:cond-capt} and \eqref{eq:cond-nsp} hold and $\rho\log\Theta<\frac14$. Then \eqref{L^2est} holds with
\[
\gamma_4:=\min\Big\{\frac{\gamma_1}{2},\ \frac{\kappa s-\log\Theta}{2},\ \frac1\rho\Big\},
\]
and $\gamma_4>\frac12\log\Theta$, so that Corollary \ref{key-cor} applies.
\end{prop}

\begin{proof}
By Lemmas \ref{lem:trapped} and \ref{lem:untrapped},
\[
\|\cL^N_\nu h_\le\|^2_{L^2}\le C\sum_\alpha\Big(\Big\|\sum_\frh a_{\alpha,\frh}\Big\|^2_{L^2}+\Big\|\sum_\frh b_{\alpha,\frh}\Big\|^2_{L^2}\Big)\le C\big(e^{-\gamma_1N}+|\nu|^{-4}\big)\|h\|^2_{L^2}.
\]
Together with \eqref{eq:tail}, and since $|\nu|^{-1}\le e^{-(N-1)/\rho}$ and $\Theta^{N/2}|\nu|^{-\varepsilon s}=\Theta^{N/2}|\nu|^{-\kappa\rho s/2}\le Ce^{-(\kappa s-\log\Theta)N/2}$,
\[
\|\cL^N_\nu h\|_{L^2}\le C\big(e^{-\gamma_1N/2}+|\nu|^{-2}+\Theta^{N/2}|\nu|^{-\varepsilon s}\big)\|h\|_{L^2}\le C_1e^{-\gamma_4N}\|h\|_{L^2}.
\]
Finally, $\gamma_1/2>\frac12\log\Theta$ and $(\kappa s-\log\Theta)/2>\frac12\log\Theta$ by \textbf{(P)}, and $1/\rho>4\log\Theta$, so $\gamma_4>\frac12\log\Theta$.
\end{proof}

By Corollary \ref{key-cor} and Proposition \ref{prop:L2est}, Theorem \ref{thm:DI} is proved.

\section{Proof of Theorem \ref{teoA}}\label{sec:proof-thmA}

We work with the complex Sobolev space $H^s$. Closed subspaces associated with the positivity sets of the invariant densities will provide the componentwise perturbation argument. Their construction does not require multiplication by characteristic functions to be bounded on $H^s$, or a full-measure union of the supports.

\subsection{Small-frequency spectral decomposition} The following lemma implements the Nagaev--Guivarc'h
spectral method in our possibly non-ergodic setting.
We include a proof sketch, explaining in particular
the decomposition associated with the ergodic acips.

\begin{lem}[Component branches and variances]\label{lem:branches}
Under the hypotheses of Theorem~\ref{teoA}, there are $\eta>0$, $C>0$, $\theta\in(0,1)$, analytic functions $\lambda_k$ and analytic rank-one projections $\Pi_k$, defined for complex $z$ near $0$, such that
\begin{equation}\label{eq:EigenSplitExact}
\cL_z^n=\sum_{k=1}^{\ell}\lambda_k(z)^n\Pi_k(z)+N(z)^n,
\qquad \|N(z)^n\|_{H^s\to H^s}\le C\theta^n.
\end{equation}
Here $\cL_z=\cL_F(e^{iz\tau}\,\cdot\,)$, $\lambda_k(0)=1$, $\Pi_k(0)u=\rho_k\ell_k(u)$, and
\[
\lambda_k'(0)=0,\qquad \lambda_k''(0)=-\sigma_k^2,
\]
where
\begin{equation}\label{eq:variance}
\sigma_k^2=\lim_{n\to\infty}\frac1n\int\tau_n^2\,\dd\mu_k
=\int\tau^2\,\dd\mu_k+2\sum_{j=1}^{\infty}\int\tau(\tau\circ F^j)\,\dd\mu_k>0.
\end{equation}
\end{lem}

\begin{proof}
Multiplication by the smooth function $\tau$ defines a bounded operator $M_\tau$ on $H^s$. Thus $\cL_z=\cL_Fe^{izM_\tau}$ is entire in operator norm. By Lemma~\ref{lem:spectral}(a), $1$ is an isolated semisimple eigenvalue of multiplicity $\ell$, and all remaining spectrum lies strictly inside the unit disk. Let $\Gamma$ be a positively oriented circle enclosing $1$
and no other spectrum of $\cL_0$. For small complex $z$,
the spectral projection
\[
P(z)=\frac{1}{2\pi i}\int_\Gamma
(w\Id-\cL_z)^{-1}\,\dd w
\]
is analytic and has rank $\ell$. To obtain the individual
analytic eigenvalue branches, we restrict $\cL_z$ to
invariant subspaces on which $1$ is a simple eigenvalue. Put $E_k=\{\rho_k>0\}$ modulo null sets and
\[
\cB_k=\{u\in H^s:u=0\text{ a.e. on }\bT^2\setminus E_k\}.
\]
This is a closed subspace of $H^s$. Positivity and $\cL_F\rho_k=\rho_k$ imply that $\cL_F$, and hence every $\cL_z$, preserves $\cB_k$. On $\cB_k$ the eigenvalue $1$ is simple: the full eigenspace at $1$ is spanned by the densities with disjoint positivity sets, and only $\rho_k$ belongs to $\cB_k$. Moreover, $\Pi|_{\cB_k}u=\rho_k\int u\,\dd\frm$ and the complementary powers decay exponentially by the full-space estimate. Analytic perturbation on this closed subspace therefore gives a simple analytic branch $\lambda_k(z)$ with eigenvector $v_k(z)\in\cB_k$, normalized by $\int v_k(z)\,\dd\frm=1$.

The $v_k(z)$ have disjoint supports modulo null sets, so are linearly independent. They belong to $\operatorname{ran}P(z)$ and form a basis there. Explicitly,
\[
\Pi_k(z)u=v_k(z)\int_{E_k}P(z)u\,\dd\frm.
\]
{The formula is valid because $\int_{E_j}v_k(z)\,\dd\frm=\delta_{jk}$.} Integration over $E_k$ is continuous on $H^s\subset L^2$. Thus these projections are analytic, they mutually vanish, commute with $\cL_z$, and sum to $P(z)$. Setting $N(z)=\cL_z(\Id-P(z))$ yields \eqref{eq:EigenSplitExact}, with uniform exponential bounds on a smaller complex neighbourhood by the resolvent formula.

For real $z$, the $L^1$ contraction of $\cL_z$ gives $|\lambda_k(z)|\le1$. Differentiating at zero and pairing with $\ell_k$ gives
\[
\lambda_k'(0)=i\ell_k\big(\cL_F(\tau\rho_k)\big)
=i\int\tau\rho_k\,\dd\frm=0.
\]
Smooth multiplication preserves $H^s$, so the spectral gap gives exponential decay of the covariances in \eqref{eq:variance}. The identity
\[
\int\tau_n^2\,\dd\mu_k
=n\int\tau^2\,\dd\mu_k+2\sum_{j=1}^{n-1}(n-j)\int\tau(\tau\circ F^j)\,\dd\mu_k
\]
proves existence of the limit. Differentiating
$\int\cL_z^n\rho_k\,\dd\frm=\lambda_k(z)^n\int\Pi_k(z)\rho_k\,\dd\frm+O(\theta^n)$ twice and dividing by $n$ identifies it with $-\lambda_k''(0)$; differentiated residual terms also decay exponentially, by analyticity and Cauchy estimates.

If $\sigma_k^2=0$, exponential covariance decay and the preceding identity give $\sup_n\|\tau_n\|_{L^2(\mu_k)}<\infty$. Write $Uu=u\circ F$ on $L^2(\mu_k)$ and set $b_M=M^{-1}\sum_{n=1}^M\tau_n$. The sequence $b_M$ is bounded, while
\[
(\Id-U)b_M=\tau-\frac1M\sum_{n=1}^MU^n\tau\longrightarrow\tau
\quad\text{in }L^2(\mu_k)
\]
by the mean ergodic theorem and centring. A weakly convergent subsequence gives $\tau=b-b\circ F$ with $b\in L^2(\mu_k)$. Consequently $e^{i\nu\tau}=\phi\circ F/\phi$ for $\phi=e^{-i\nu b}$ and every real $\nu\ne0$, contradicting \textbf{(B2)}. This proves strict positivity.
\end{proof}

\begin{lem}[Characteristic-function expansion]\label{lem:Char}
For $\psi\in H^s$, $\xi\in L^2$, every integer $R\ge0$, and every fixed $b>0$, there are polynomials $A_{j,k,\psi,\xi}$, independent of $n$, with degree at most $3j$ and parity $j$, such that
\begin{equation}\label{CharAsympExp}
\int\psi e^{it\tau_n/\sqrt n}\xi\circ F^n\,\dd\frm
=\sum_{k=1}^{\ell}e^{-\sigma_k^2t^2/2}
\sum_{j=0}^{R}n^{-j/2}A_{j,k,\psi,\xi}(t)+E_{R,n}(t)
\end{equation}
for $|t|\le b\sqrt{\log n}$. Here $A_{0,k,\psi,\xi}=c_k(\psi,\xi)$, and, for some $c>0$, integer $M_R$, and $\theta_1\in(0,1)$,
\begin{equation}\label{CharError}
|E_{R,n}(t)|\le C n^{-(R+1)/2}(1+|t|^{M_R})e^{-ct^2}+C\theta_1^n.
\end{equation}
For $\psi=v$ a probability density and $\xi=1$, one has $A_{j,k,v,1}(0)=0$ for $j\ge1$, and the improved estimate
\begin{equation}\label{eq:char-error-normalized}
|E_{R,n}(t)|\le C n^{-(R+1)/2}|t|(1+|t|^{M_R})e^{-ct^2}
+C\theta_1^n|t|.
\end{equation}
For sufficiently small fixed $\eta>0$, the characteristic function in \eqref{CharAsympExp} is bounded by $C(e^{-ct^2}+\theta_1^n)$ on $|t|\le\eta\sqrt n$.
\end{lem}

\begin{proof}
Set $a_k(z)=\int\xi\Pi_k(z)\psi\,\dd\frm$. The scalar spectral decomposition is
\[
\int\psi e^{iz\tau_n}\xi\circ F^n\,\dd\frm
=\sum_k\lambda_k(z)^na_k(z)+\int\xi N(z)^n\psi\,\dd\frm.
\]
The analytic logarithm near $\lambda_k(0)=1$ has the expansion
\[
\log\lambda_k(z)=-\tfrac12\sigma_k^2z^2+\sum_{r=3}^{R+2}\beta_{r,k}z^r+O(|z|^{R+3}).
\]
Insert $z=t/\sqrt n$, expand $a_k$ and the exponential, and collect powers of $n^{-1/2}$. A log-eigenvalue term of degree $r$ has order $n^{-(r-2)/2}$ and an amplitude term of degree $r$ has order $n^{-r/2}$. This proves the degree bound and parity. Taylor remainders on $|t|\le b\sqrt{\log n}$ are bounded by a fixed polynomial times $n^{-(R+1)/2}e^{-ct^2}$. The spectral residual is exponentially small. This proves \eqref{CharError}.

For $\xi=1$ and probability density $v$, the approximating characteristic function equals $1$ at $t=0$: each $a_k(0)=c_k(v)$, $\sum_kc_k(v)=1$, and every correction polynomial vanishes at zero. The residual integral also vanishes at zero, since $\int N(0)^nv\,\dd\frm=0$. Applying Taylor's formula to these remainders, with differentiated analytic estimates, inserts a factor $|t|$ and proves \eqref{eq:char-error-normalized}. Finally, positivity of the variances gives $|\lambda_k(z)|\le e^{-cz^2}$ for real $|z|\le\eta$, proving the last assertion.
\end{proof}

\subsection{Intermediate and large frequencies}

\begin{lem}[A uniform scalar bound]\label{lem:uniform-scalar}
There is $C>0$ such that for all real $\nu$, all $n\ge0$, $\psi\in H^s$, and $\xi\in L^2$,
\begin{equation}\label{eq:uniform-scalar}
\int |\psi|\,|\xi\circ F^n|\,\dd\frm
\le C\|\psi\|_{H^s}\|\xi\|_{L^2},
\qquad
\left|\int\xi\cL_\nu^n\psi\,\dd\frm\right|
\le C\|\psi\|_{H^s}\|\xi\|_{L^2}.
\end{equation}
\end{lem}

{\begin{proof}
If $0<s\le1$, the modulus map satisfies $\||\psi|\|_{H^s}\le C\|\psi\|_{H^s}$, by the fractional difference norm or the weak chain rule. Power boundedness of $\cL_F$ on $H^s$ therefore gives
$\|\cL_F^n|\psi|\|_2\le C\|\psi\|_{H^s}$.
If $s>1$, Sobolev embedding in dimension two gives instead
\[
\|\cL_F^n|\psi|\|_2\le\|\psi\|_\infty\|\cL_F^n1\|_2
\le C\|\psi\|_{H^s},
\]
again by \eqref{LY-F}. In either case, positivity and transfer-operator duality yield
\[
\int |\psi|\,|\xi\circ F^n|\,\dd\frm
=\int |\xi|\cL_F^n|\psi|\,\dd\frm
\le C\|\psi\|_{H^s}\|\xi\|_2.
\]
For unbounded $\xi$ the equality follows by monotone truncation. Finally $|\cL_\nu^n\psi|\le\cL_F^n|\psi|$ proves the second bound.
\end{proof}}

For every fixed $K>0$, $B>0$, and integer $R\ge1$, Lemma~\ref{lem:spectral}(c) and Theorem~\ref{thm:DI} imply
\begin{equation}\label{eq:C4}
\int_{K<|\nu|<Bn^{(R-1)/2}}
\frac{\|\cL_\nu^n\|_{H^s\to H^s}}{|\nu|}\,\dd\nu
=o(n^{-R/2}).
\end{equation}
Indeed, the compact part bounded away from zero is exponentially small, while the remaining integral is at most
\[
Ce^{-\gamma n}\int_{\nu_0}^{\max\{\nu_0,Bn^{(R-1)/2}\}}
\nu^{s+\delta-1}\,\dd\nu,
\]
which is an exponential times a fixed power of $n$. Note that the restriction away from zero is essential here, since a uniform exponential spectral bound on $0<|\nu|\le\nu_0$ would contradict $\lambda_k(\nu)\to1$.

For later use, if $p>1$ and $g\in\fF_0^p$, the integral over $|\nu|>n^a$ in Fourier inversion is at most
\begin{equation}\label{eq:far-tail}
C\|\psi\|_{H^s}\|\xi\|_2 C^p(g)n^{-a(p-1)}
\end{equation}
by Lemma~\ref{lem:uniform-scalar}. On $\eta\le|\nu|\le n^a$, the spectral compact-frequency bound and DI make the integral smaller than any inverse power of $n$, for each fixed $a$. Taking $a$ large proves that every fixed $q>0$ in the test-function classes suffices for every fixed expansion order.

\subsection{Edgeworth expansions}

Let $v\in H^s$ be nonnegative with $\int v\,\dd\frm=1$, and let $\chi_n(t)=\int v e^{it\tau_n/\sqrt n}\,\dd\frm$. Define
\begin{equation}\label{eq:kernel-definition}
K_{j,k,\psi,\xi}(x)=\frac1{2\pi}\int_{\bR}e^{-itx}
e^{-\sigma_k^2t^2/2}A_{j,k,\psi,\xi}(t)\,\dd t,
\qquad K_{j,\psi,\xi}=\sum_kK_{j,k,\psi,\xi}.
\end{equation}
Each kernel is a polynomial of degree at most $3j$ times $\fn_{\sigma_k}$; this formula fixes the Fourier signs and normalizations. In the case $(\psi,\xi)=(v,1)$, $\int K_{j,k,v,1}=0$ for $j\ge1$. Hence there is a polynomial $P_{j,k,v}$ satisfying
\[
\frac{\dd}{\dd x}\big(P_{j,k,v}(x)\fn_{\sigma_k}(x)\big)=K_{j,k,v,1}(x).
\]
Define $\cE_{R,n}$ by the two explicit sums on the right-hand side of \eqref{GenEdgeExp}. Its derivative is a signed density of total mass one, uniformly bounded in $n$, and its Fourier--Stieltjes transform with the positive characteristic-function sign is
\[
\chi_{R,n}(t)=\sum_k e^{-\sigma_k^2t^2/2}
\left(c_k(v)+\sum_{j=1}^R n^{-j/2}A_{j,k,v,1}(t)\right).
\]
The Gaussian factor multiplies every correction term. The usual smoothing inequality for comparison with such a signed distribution function gives
\[
\sup_z|\bP_{v\frm}(\tau_n/\sqrt n\le z)-\cE_{R,n}(z)|
\le C\int_{-T}^T\frac{|\chi_n(t)-\chi_{R,n}(t)|}{|t|}\,\dd t+\frac{C}{T};
\]
see the Edgeworth smoothing argument in \cite{FeLi,FePe}.

Choose $T=n^{(R+1)/2}$. On $|t|\le b\sqrt{\log n}$, \eqref{eq:char-error-normalized} makes the integral $O(n^{-(R+1)/2})+O(\theta_1^n\sqrt{\log n})$. On $b\sqrt{\log n}<|t|\le\eta\sqrt n$, the Gaussian estimates from Lemma~\ref{lem:Char}, and the same estimates for the approximating polynomials, make the integral $o(n^{-R/2})$ when $b$ is large enough. On $\eta\sqrt n<|t|\le T$, substitute $\nu=t/\sqrt n$. The actual characteristic function is controlled by \eqref{eq:C4}, with its polynomial upper cutoff $n^{R/2}$ (equivalently, use that estimate with order $R+1$). The approximating characteristic function has an exponentially small Gaussian tail. Finally $T^{-1}=o(n^{-R/2})$. This proves \eqref{GenEdgeExp} uniformly in $z$.

\subsection{Mixing global and local expansions}

Let $p=q+2>2$. Since $\widehat g\in L^1$ and Lemma~\ref{lem:uniform-scalar} bounds the absolute endpoint integral, Fubini's theorem and Fourier inversion give
\begin{equation}\label{eq:fourier-mixing}
I_n(g;\psi,\xi)=\frac1{2\pi}\int_{\bR}\widehat g(\nu)
\left(\int\xi\cL_\nu^n\psi\,\dd\frm\right)\,\dd\nu.
\end{equation}
By \eqref{eq:far-tail}, the compact-frequency bound and DI, all frequencies outside $|\nu|\le\eta$ contribute $o(n^{-L})$ for any prescribed $L$. After $t=\sqrt n\,\nu$, frequencies with $b\sqrt{\log n}<|t|\le\eta\sqrt n$ also contribute $o(n^{-L})$ if $b$ is sufficiently large. On the remaining interval insert \eqref{CharAsympExp}. The integrated error in \eqref{CharError}, including the factor $n^{-1/2}$ from the change of variable, is $O(n^{-(R+2)/2})+O(\theta_1^n)$. Extending the Gaussian integrals to $\bR$ gives
\[
I_n(g;\psi,\xi)=\sum_{j=0}^{R}\frac{n^{-(j+1)/2}}{2\pi}
\sum_k\int_{\bR}\widehat g(t/\sqrt n)
A_{j,k,\psi,\xi}(t)e^{-\sigma_k^2t^2/2}\,\dd t
+o(n^{-(R+1)/2}).
\]
Using \eqref{eq:kernel-definition} and Fubini, the integral divided by $2\pi$ equals
$\int g(x)K_{j,\psi,\xi}(x/\sqrt n)\,\dd x$.
This proves the global expansion \eqref{eq:global-expansion}. In particular, its leading size is $n^{-1/2}$, as required by the local limit theorem.

For the local expansion suppose $g\in\fF_{R+1}^p$. Multiply the global expansion by $\sqrt n$ and Taylor expand each $K_j(x/\sqrt n)$ at zero through degree $R-j$. All derivatives of these Gaussian-polynomial kernels are bounded. Thus the integrated Taylor errors are $O(n^{-(R+1)/2}M_{R+1}(g))$. The parity in Lemma~\ref{lem:Char} implies $K_j(-x)=(-1)^jK_j(x)$; hence the terms with $j+d$ odd vanish. Define
\begin{equation}\label{eq:Q-definition}
Q_{m,\psi,\xi}(x)=\sum_{j+d=2m}\frac{K_{j,\psi,\xi}^{(d)}(0)}{d!}x^d.
\end{equation}
Collecting the even powers yields \eqref{eq:local-expansion}, with $\deg Q_m\le2m$ and the stated leading coefficient. This completes the proof of Theorem~\ref{teoA}.

\section{Proof of Theorem \ref{teoB}}\label{sec:doc}

{We use the fibrewise Fourier argument of \cite{BuEs}.} Recall that $\tilde F(x,y,z)=(F(x,y), z+\tau(x,y)) \mod \bZ^3$. Thanks to the skew-product structure, we can exploit the fiber-wise Fourier expansion of a function $f\in\cC^r(\bT^3,\bC)$:
\[
f(x,y,z)=\sum_{k\in \bZ} e^{2\pi i k z} \hat{f}_k(x,y),\quad \hat{f}_k(x,y)=\int_{\bT} e^{-2\pi i k z} f(x,y,z) \, \dd z.
\]
Integrating by parts in $z$, for $0\le j\le r$ we have
\begin{equation}\label{eq:fourier-decay}
\|\hat f_k\|_{\cC^j(\bT^2)}\le \|f\|_{\cC^j},\qquad \|\hat f_k\|_{\cC^0(\bT^2)}\le (2\pi|k|)^{-r}\|f\|_{\cC^r}\quad (k\ne0).
\end{equation}
For notational convenience we set $p=(x,y)\in \bT^2$ and $\dd p= \dd x\dd y$. We will prove the statement for observables in $\cC^r$, with an integer $r \ge \max\{\lceil s\rceil,2\}$, and deduce the conclusion for Hölder observables at the end of the proof. By assumption $\mu_*$ is an ergodic acip of $F$, hence $\mu_*=\mu_{k_0}$ and $h_*=\rho_{k_0}\in H^{s}(\bT^2)$ for some $k_0$. Let $g,h\in \cC^{r}(\bT^3,\bC)$ and assume that ${\mu}(g)=0$, where $\dd{\mu}=\dd{\mu}_*\times \operatorname{Leb}$; note in particular that this implies $\mu_*(\hat g_0)=0$. Recall that, for $u\in L^1$, $w\in L^\infty$ and $\nu\in\bR$,
\[
\int_{\bT^2} u\cdot (w\circ F^n)\, e^{i\nu\tau_n}\,\dd p=\int_{\bT^2}\cL^n_\nu(u)\, w\,\dd p.
\]
We compute for each $n\in \bN$, by linearity and absolute uniform convergence of the Fourier series (as $r\ge2$),
\[
\begin{split}
    &\int_{\bT^3} (g\cdot h\circ \tilde F^n)(p,z)h_*(p)\,\dd p\,\dd z\\
    &\quad =\sum_{k\in \bZ} \sum_{j\in \bZ} \int_{\bT^2}(\hat g_{k}\cdot h_*)(p)\, \hat h_{j}(F^n(p))\,e^{2\pi i j \tau_n(p)}\dd p \int_{\bT}e^{2\pi i(k+j) z}\, \dd z\\
    &\quad =\sum_{k\in \bZ}\int_{\bT^2}  \cL_{-2\pi k}^n(\hat g_{k}\cdot h_*)(p)\,\hat h_{-k}(p)\, \dd p,
\end{split}
\]
where in the last line we used that $\int_{\bT}e^{2\pi i(k+j)z}\dd z=0$ unless $j=-k$.

Set $A_k:=\int_{\bT^2} \cL_{-2\pi k}^n(\hat g_{k}\cdot h_*)(p)\,\hat h_{-k}(p)\, \dd p$ and write
\[
\sum_{k\in \bZ} A_k=A_0+\sum_{0<2\pi|k|< \nu_0} A_k+ \sum_{2\pi|k|\ge\nu_0} A_k,
\]
{where $\nu_0$ is given by Theorem \ref{thm:DI}, enlarged so that $\nu_0>2\pi$.} We prove that the absolute value of each of these three terms decays exponentially. Below we use Lemma \ref{lem:leib} in the form $\|\varphi\, u\|_{H^s}\le C\|\varphi\|_{\cC^{\lceil s\rceil}}\|u\|_{H^s}$.

\emph{The term $A_0$.} Let $\Pi$ be the spectral projector of $\cL_F$ on $H^s$ associated with the eigenvalue $1$. By Lemma \ref{lem:spectral}(a) there exist $C>0$, $\theta_0\in(0,1)$ such that $\|\cL_F^n(\Id-\Pi)\|_{H^s}\le C\theta_0^n$. We claim that $\Pi(\hat g_0h_*)=0$. Let $E:=\{\rho_{k_0}>0\}$. Since $\rho_{k_0}\ge0$ and $\cL_F\rho_{k_0}=\rho_{k_0}$, {for almost every $x\notin E$ we have $\rho_{k_0}(y)=0$ for all $y\in F^{-1}(x)$;} hence, if $u$ vanishes outside $E$, so does $\cL_Fu$. By Lemma \ref{lem:spectral}(a), $\Pi u=\lim_n\cL_F^nu$ in $H^s$, so $\Pi u$ vanishes outside $E$ as well. Since $\Pi u$ is a combination of the $\rho_j$,{whose positivity sets are disjoint modulo null sets}, $\Pi u=c\,\rho_{k_0}$ with $c=\int\Pi u=\lim_n\int\cL_F^nu=\int u$. Applying this to $u=\hat g_0h_*=\hat g_0\rho_{k_0}$, for which $\int u=\mu_*(\hat g_0)=0$, proves the claim. Hence
\[
\begin{split}
    |A_0|
&=\left|\int_{\bT^2} \cL_F^n(\hat g_0 h_*)(p)\,\hat h_0(p)\,\dd p\right|
\le \|\cL_F^n(\Id-\Pi)(\hat g_0 h_*)\|_{L^2}\,\|\hat h_0\|_{L^2}\\
&\le C\,\theta_0^n \|\hat g_0 h_*\|_{H^{s}}\|\hat h_0\|_{L^2}\le  C\,\theta_0^n \| g\|_{\cC^r}\|h\|_{\cC^r}\|h_*\|_{H^s}.
\end{split}
\]

\emph{The terms $0<2\pi|k|<\nu_0$.} These are finitely many, and by \eqref{midOpEst} with $J=\{2\pi\le|\nu|\le\nu_0\}$,
\[
\max_{0<2\pi|k|< \nu_0}|A_k| \le C \theta_J^n \max_{0<2\pi|k|< \nu_0} \|\hat g_k h_*\|_{H^s}\|\hat h_{-k}\|_{L^2} \leq C\theta_J^n\| g\|_{\cC^r}\|h\|_{\cC^r}\|h_*\|_{H^s}.
\]

\emph{The terms $2\pi|k|\ge\nu_0$.} By \eqref{main-est}, $\|u\|_{L^2}\le\|u\|_\nu\le 2\|u\|_{H^s}$, and \eqref{eq:fourier-decay},
\[
\begin{split}
|A_k|&\le \|\cL^n_{-2\pi k}(\hat g_k h_*)\|_{-2\pi k}\,\|\hat h_{-k}\|_{L^2}\le C|k|^{\delta}e^{-\gamma n}\|\hat g_k h_*\|_{H^s}\,\|\hat h_{-k}\|_{L^2}\\
&\le C\,|k|^{\delta-r}\,e^{-\gamma n}\,\|g\|_{\cC^r}\|h\|_{\cC^r}\|h_*\|_{H^s},
\end{split}
\]
and therefore
\[
\bigg|\sum_{2\pi|k|\ge\nu_0} A_k \bigg| \le C e^{-\gamma n}\|g\|_{\cC^r}\|h\|_{\cC^r}\|h_*\|_{H^s} \sum_{k\ne0} |k|^{\delta-r},
\]
where the last sum converges since $r\ge2>1+\delta$. This proves the statement for $\cC^r$ observables, with some rate $\gamma'>0$ and constant $C\|g\|_{\cC^r}\|h\|_{\cC^r}$.

\emph{Hölder observables.} Let $g,h\in\cC^\alpha(\bT^3)$, $\alpha\in(0,1)$, and for $\epsilon\in(0,1)$ let $g_\epsilon,h_\epsilon$ be their convolutions with a smooth mollifier at scale $\epsilon$, so that
\[
\|g-g_\epsilon\|_{\cC^0}\le \epsilon^\alpha\|g\|_{\cC^\alpha},\qquad \|g_\epsilon\|_{\cC^r}\le C\epsilon^{\alpha-r}\|g\|_{\cC^\alpha},
\]
and similarly for $h$. Since $\mu$ is $\tilde F$-invariant, replacing $(g,h)$ by $(g_\epsilon,h_\epsilon)$ changes both $\mu(g\cdot h\circ\tilde F^n)$ and $\mu(g)\mu(h)$ by at most $C\epsilon^\alpha\|g\|_{\cC^\alpha}\|h\|_{\cC^\alpha}$. Applying the $\cC^r$ estimate to $g_\epsilon-\mu(g_\epsilon)$ and $h_\epsilon$,
\[
\big|\mu(g\cdot h\circ\tilde F^n)-\mu(g)\mu(h)\big|\le C\big(\epsilon^\alpha+\epsilon^{2(\alpha-r)}e^{-\gamma'n}\big)\|g\|_{\cC^\alpha}\|h\|_{\cC^\alpha},
\]
and the choice $\epsilon=e^{-\gamma'n/(2r)}$ yields the bound $Ce^{-\alpha\gamma'n/(2r)}\|g\|_{\cC^\alpha}\|h\|_{\cC^\alpha}$. \qed

For $\alpha\ge1$, the assertion follows by the continuous embedding into $\cC^\beta$ for any fixed $0<\beta<1$.
\appendix

\section{Partial captivity in the case of skew products}\label{app:B1}

In this Appendix we characterise \textbf{(B1)} for skew products and provide an example to which our results apply. Let
\[
F(x,\theta)=(f(x),\theta+\omega(x)),\qquad \tau=\tau(x),
\]
where $f\in\cC^\infty(\bT,\bT)$ is expanding of degree $d\ge2$, $\omega,\tau\in\cC^\infty(\bT,\bR)$, and $F$ is an SVPH. Let $h\in\cC^\infty(\bT)$ be the positive invariant density of $f$ ($\cL_fh=h$, $\int h=1$). For $x\in\bT$ and $n\in\bN$ let $\cH_f^n(x)$ be the set of inverse branches $g$ of $f^n$ near $x$. Writing $\Phi:=(\omega,\tau):\bT\to\bR^2$ and $\Phi_n:=\sum_{k=0}^{n-1}\Phi\circ f^k$,\footnote{Similarly we write
$
\omega_n=\sum_{j=0}^{n-1}\omega\circ f^j$ and $
\tau_n=\sum_{j=0}^{n-1}\tau\circ f^j.
$} we set, for $g\in\cH^n_f(x)$,
\[
{\begin{gathered}P_g(x):=(\Phi_n\circ g)'(x)=\big((\omega_n\circ g)'(x),(\tau_n\circ g)'(x)\big)=:(s_g,t_g)\in\bR^2,\\ p_g(x):=\frac{g'(x)\,h(g(x))}{h(x)},\end{gathered}}
\]
and $\mu^x_n:=\sum_{g\in\cH^n_f(x)}p_g(x)\,\delta_{P_g(x)}$, a probability measure on $\bR^2$. We say that a function $\phi\in\cC^\infty(\bT)$ is \emph{cohomologous to a constant} if $\phi=C\circ f-C+k$ for some $C\in\cC^1(\bT)$ and $k\in\bR$.{\footnote{By measurable Liv\v{s}ic regularity for expanding maps, a measurable solution (with the identity holding almost everywhere for the smooth invariant probability) has a smooth representative; in particular the transfer function may be taken in $\cC^\infty$.}}

\begin{lem}\label{lem:skew}
Let $F,\tau$ be as above and assume that $\omega$ is not cohomologous to a constant.\footnote{This is necessary for $F$ to be VE: if $\omega=\psi\circ f-\psi+c$, then $F$ is smoothly conjugate to $(x,\theta)\mapsto(f(x),\theta+c)$, for which $B_s(F^n)\ge\inf\cL_f^n1$ (take $\xi=d\theta$ in \eqref{Bs}), so that \eqref{virtually} fails.} Then the following are equivalent:
\begin{enumerate}[label=(\roman*)]
\item \textbf{(B1)} holds for some $\kappa\in(0,\log\lambda)$;
{\item \textbf{(B1)} holds for every $\kappa\in(0,\log\lambda)$, with $C$ and $\gamma_1$ allowed to depend on $\kappa$;}
\item for every $b\in\bR$, $\tau-b\,\omega$ is not cohomologous to a constant.
\end{enumerate}
Moreover, if these hold, then $F$ satisfies \eqref{main-hp}.
\end{lem}

Condition (iii) with $b=0$ is the non-coboundary condition on the fibre map $\tau$. For $b\neq0$ it excludes a resonance between the two fibre cocycles $\omega$ and $\tau$. Together with the assumption on $\omega$, (iii) says that $a\omega+b\tau$ is not cohomologous to a constant for every $(a,b)\neq(0,0)$, which is the condition under which exponential mixing of the $\bT^2$-extension $\tilde F$ of $f$ is proved in \cite{ChHu}. Finally, the assumption that $\omega$ is not cohomologous to a constant is the condition that makes the unstable cones of $F$ transversal, see \eqref{main-hp}.

\begin{proof}
Throughout, $C$ denotes constants depending only on $f,\omega,\tau$. By \eqref{partial hyperbolicity 2} applied to $v=(1,0)$, $C^{-1}\Lambda^{-n}\le g'(x)\le C\lambda^{-n}$ for $g\in\cH^n_f(x)$, and $C^{-1}g'(x)\le p_g(x)\le Cg'(x)$.

\emph{Step 1: translation of \eqref{eq:B1star}.} Let $q=(x,\theta)$. Since $F^n(y,\phi)=(f^n(y),\phi+\omega_n(y))$, the inverse branches of $F^n$ at $q$ are $\frh(x,\theta)=(g(x),\theta-\omega_n(g(x)))$, $g\in\cH^n_f(x)$, and $J_\frh(q)=g'(x)$, $\tau_n\circ\frh(x,\theta)=\tau_n(g(x))$. With $y=g(x)$,
\[
D_{\frh(q)}F^n\begin{pmatrix}1\\ \sigma\end{pmatrix}=\begin{pmatrix}(f^n)'(y)\\ \omega_n'(y)+\sigma\end{pmatrix},
\]
{so, up to sign (which does not affect the absolute-value condition), the unit vectors of $D_{\frh(q)}F^n\fC_u$ are $u=(1,s)/\sqrt{1+s^2}$} with $s=s_g+\sigma g'(x)$, $|\sigma|\le\chi_u$; in particular $|s|\le\chi_u<1$. Moreover $d(\tau_n\circ\frh)_q(u)=t_g/\sqrt{1+s^2}$ and $\zeta(u)=(\zeta_1+\zeta_2s)/\sqrt{1+s^2}$ for $\zeta=(\zeta_1,\zeta_2)$. Let $\delta_n:=e^{\kappa n}\lambda^{-n}$ and $\ell_\zeta:=\{(s,t):t=\zeta_1+\zeta_2s\}$. Since $\operatorname{dist}(P_g,\ell_\zeta)=|\zeta_1+\zeta_2s_g-t_g|/\sqrt{1+\zeta_2^2}$, we get:
\begin{itemize}
\item if $\frh$ belongs to the set in \eqref{eq:B1star}, then $|\zeta_1+\zeta_2s_g-t_g|\le\sqrt2\,\delta_n+|\zeta_2|\chi_ug'(x)$, hence $\operatorname{dist}(P_g,\ell_\zeta)\le\sqrt2\,\delta_n+C\lambda^{-n}\le C\delta_n$;
\item if $|\zeta_1+\zeta_2s_g-t_g|\le\delta_n$, then $\frh$ belongs to the set in \eqref{eq:B1star} (take $\sigma=0$).
\end{itemize}
Consequently, with $\ell^{(r)}$ the $r$-neighbourhood of $\ell$,
\begin{equation}\label{eq:skew-transl}
\sum_{g:\ |\zeta_1+\zeta_2s_g-t_g|\le\delta_n}g'(x)\ \le\ \text{(left side of \eqref{eq:B1star})}\ \le\ C\,\mu^x_n\big(\ell_\zeta^{(C\delta_n)}\big).
\end{equation}

\emph{Step 2: (i)$\Rightarrow$(iii).} Suppose $\tau-b\omega=C_0\circ f-C_0+k$ with $C_0\in\cC^1$. Then $\tau_n-b\omega_n=C_0\circ f^n-C_0+nk$, so $t_g-bs_g=\big((\tau_n-b\omega_n)\circ g\big)'(x)=C_0'(x)-C_0'(g(x))g'(x)$ and $|C_0'(x)+bs_g-t_g|\le C\lambda^{-n}\le\delta_n$ for $n$ large. Taking $\zeta=(C_0'(x),b)$ in \eqref{eq:skew-transl}, the left side of \eqref{eq:B1star} is at least $\sum_gg'(x)=\cL_f^n1(x)\ge C^{-1}$, so \textbf{(B1)} fails for every $\kappa>0$.

\emph{Step 3: the limit measures.} For $g\in\cH^n_f(x)$ and $\tilde g\in\cH^m_f(g(x))$ we have $\tilde g\circ g\in\cH^{n+m}_f(x)$ and, since $\Phi_{n+m}\circ\tilde g\circ g=\Phi_m\circ\tilde g\circ g+\Phi_n\circ g$,
\begin{equation}\label{eq:skew-ss}
P_{\tilde g\circ g}(x)=P_g(x)+g'(x)\,P_{\tilde g}(g(x)),\qquad p_{\tilde g\circ g}(x)=p_g(x)\,p_{\tilde g}(g(x)).
\end{equation}
Moreover $|P_g(x)|\le\sum_{k=1}^n\|\Phi'\|_\infty\sup|(f^{k})'|^{-1}\le C_0$. Hence $$\mu^x_{n+m}=\sum_{g\in\cH^n_f(x)}p_g(x)\,(S_g)_*\mu^{g(x)}_m,$$ where $S_g(P):=P_g(x)+g'(x)P$, and each $(S_g)_*\mu_m^{g(x)}$ is supported in the ball of radius $C_0g'(x)\le C\lambda^{-n}$ around $P_g(x)$. Thus $(\mu_n^x)_n$ is Cauchy in the Wasserstein distance, uniformly in $x$, and converges to a probability measure $\mu^x$, depending continuously on $x$, with
\begin{equation}\label{eq:skew-ss2}
\mu^x=\sum_{g\in\cH^n_f(x)}p_g(x)\,(S_g)_*\mu^{g(x)},\qquad \mu^x_n(A)\le\mu^x\big(A^{(C\lambda^{-n})}\big){\quad\text{for all }n\text{ and Borel }A\subset\bR^2} .
\end{equation}

\emph{Step 4: no line carries $\mu^x$.} We claim that, under (iii) and the assumption on $\omega$, there is no $x\in\bT$ and no line $\ell\subset\bR^2$ with $\mu^x(\ell)=1$. Suppose $\operatorname{supp}\mu^{x_0}\subset\{\alpha s+\beta t=c_0\}$, with $\alpha^2+\beta^2=1$. For $y\in\bT$ let $c(y):=\int(\alpha s+\beta t)\,\dd\mu^y$ and $V(y):=\int(\alpha s+\beta t-c(y))^2\dd\mu^y$; both are continuous and $V(x_0)=0$. By \eqref{eq:skew-ss2} with $n=1$, and since $p_g>0$ and $S_g$ is a homothety, $V(x_0)=0$ implies $V(z)=0$ for every $z\in f^{-1}(x_0)$; by induction $V=0$ on $\bigcup_nf^{-n}(x_0)$, which is dense, hence $V\equiv0$: for every $y$, $\operatorname{supp}\mu^y\subset\{\alpha s+\beta t=c(y)\}$. Now let $z\in\bT$, $y=f(z)$, and let $g\in\cH^1_f(y)$ with $g(y)=z$, so that $P_g(y)=\Phi'(z)/f'(z)$ and $g'(y)=1/f'(z)$. For $P=(s,t)\in\operatorname{supp}\mu^z$ we have $S_gP\in\operatorname{supp}\mu^y$, i.e.
\[
\frac{\alpha(\omega'(z)+s)+\beta(\tau'(z)+t)}{f'(z)}=c(f(z)),\qquad \alpha s+\beta t=c(z),
\]
so that $\alpha\omega'+\beta\tau'=(c\circ f)\,f'-c$ on $\bT$. Let $\tilde f$ be a lift of $f$ and $\tilde C$ a primitive of $c$ on $\bR$; then $\alpha\omega+\beta\tau=\tilde C\circ\tilde f-\tilde C+k$ for some $k\in\bR$. Since $\tilde C(x+1)=\tilde C(x)+\int_\bT c$ and $\tilde f(x+1)=\tilde f(x)+d$, periodicity of the left side forces $(d-1)\int_\bT c=0$, i.e.\ $\tilde C$ is $1$-periodic. Hence $\alpha\omega+\beta\tau$ is cohomologous to a constant. If $\beta=0$ this contradicts the assumption on $\omega$; if $\beta\ne0$ it contradicts (iii) with $b=-\alpha/\beta$.

\emph{Step 5.} Let $\vartheta(\delta):=\sup_{x\in\bT}\sup_\ell\mu^x(\ell^{(\delta)})$, the supremum over all lines $\ell\subset\bR^2$. We show that, under (iii), there are $C,\beta_0>0$ with
\begin{equation}\label{eq:skew-nc}
\vartheta(\delta)\le C\delta^{\beta_0}\qquad\forall\,\delta\in(0,1].
\end{equation}
First, there exist $N\in\bN$ and $\varepsilon>0$ such that $\mu^x(\ell^{(3C\lambda^{-N})})\le1-\varepsilon$ for all $x$ and $\ell$. Otherwise there are $N_k\to\infty$, $x_k\in\bT$ and lines $\ell_k$ with $\mu^{x_k}(\ell_k^{(3C\lambda^{-N_k})})\ge1-\frac1k$; since all $\mu^x$ are supported in $\{|P|\le C_0\}$, the lines $\ell_k$ meet a fixed ball, and along a subsequence $x_k\to x$, $\ell_k\to\ell$. By weak continuity, $\mu^x(\ell^{(\eta)})=1$ for every $\eta>0$, i.e.\ $\mu^x(\ell)=1$, contradicting Step 4.
Next, fix $\delta\le C\lambda^{-N}$, $x$ and $\ell$. By \eqref{eq:skew-ss2} with $n=N$,
\[
\mu^x(\ell^{(\delta)})=\sum_{g\in\cH^N_f(x)}p_g(x)\,\mu^{g(x)}\big(S_g^{-1}\ell^{(\delta)}\big).
\]
Since $(S_g)_*\mu^{g(x)}$ is supported in $B(P_g,C\lambda^{-N})$, only the $g$ with $\operatorname{dist}(P_g,\ell)\le2C\lambda^{-N}$ contribute; for these the ball $B(P_g,C\lambda^{-N})$ lies in $\ell^{(3C\lambda^{-N})}$, so their total weight is at most $\mu^x(\ell^{(3C\lambda^{-N})})\le1-\varepsilon$. Moreover $S_g^{-1}\ell^{(\delta)}=(\ell')^{(\delta/g'(x))}$ for a line $\ell'$ parallel to $\ell$, and $\delta/g'(x)\le C\Lambda^N\delta=:K\delta$. Hence $\vartheta(\delta)\le(1-\varepsilon)\vartheta(K\delta)$ for $\delta\le C\lambda^{-N}$, and iterating gives \eqref{eq:skew-nc} with $\beta_0=\log(1-\varepsilon)^{-1}/\log K$.

\emph{Step 6: (iii)$\Rightarrow$(ii), and \eqref{main-hp}.} Let $\kappa\in(0,\log\lambda)$. By \eqref{eq:skew-transl}, \eqref{eq:skew-ss2} and \eqref{eq:skew-nc}, the left side of \eqref{eq:B1star} is bounded by
\[
C\mu_n^x\big(\ell_\zeta^{(C\delta_n)}\big)\le C\mu^x\big(\ell_\zeta^{(C\delta_n+C\lambda^{-n})}\big)\le C\vartheta(C\delta_n)\le Ce^{-\beta_0(\log\lambda-\kappa)n},
\]
which is \textbf{(B1)} with $\gamma_1=\beta_0(\log\lambda-\kappa)$. For \eqref{main-hp}: if two branches $\frh,\tilde\frh$ at $q$ are not transversal, the slope intervals $[s_g\pm\chi_ug']$ and $[s_{\tilde g}\pm\chi_u\tilde g']$ intersect, so $|s_g-s_{\tilde g}|\le C\lambda^{-n}$; hence $\cN_F(n)\le C\sup_{x}\sup_{s_0}\mu_n^x(\{|s-s_0|\le C\lambda^{-n}\})\le C\vartheta(C\lambda^{-n})\le C\lambda^{-\beta_0n}$.

Since (ii)$\Rightarrow$(i) is trivial, the proof is complete.
\end{proof}
We conclude with a simple example satisfying all our assumptions.

\begin{ex}\label{ex:B1star}
Let $f(x)=\ell x$ with $\ell\ge2$ an integer, $\omega(x)=a\cos(2\pi x)$ with $a\neq0$, and $\tau(x)=\sin(2\pi x)$, where $\ell>2(1+2\pi|a|)$, so that \textbf{(A5)} holds, and $\ell$ is large enough for $F$ to be an SVPH.

\emph{Cohomology.} If $\phi=C\circ f-C+k$ with $C\in L^2(\bT)$, then $\hat\phi(m)=\hat C(m/\ell)\mathds 1_{\ell\mid m}-\hat C(m)$ for $m\ne0$, and therefore $\sum_{j\ge0}\hat\phi(\ell^j)=-\lim_{J\to\infty}\hat C(\ell^J)=0$. For $\phi=\tau-b\omega$ this sum equals $\hat\phi(1)=\frac1{2i}-\frac{ab}2\ne0$ for every $b\in\bR$, and for $\phi=\omega$ it equals $\frac a2\neq0$. Hence $\omega$ is not cohomologous to a constant and (iii) of Lemma \ref{lem:skew} holds, so $\tau$ satisfies \textbf{(B1)}.

\emph{Multiplicative cohomology.} Let $\varphi\in\cC^\infty(\bT,\bR)$ and suppose that $e^{i\varphi}=c\,\psi\circ f/\psi$ Lebesgue-a.e., with $c\in S^1$ and $\psi:\bT\to S^1$ measurable. By Livšic regularity for expanding maps $\psi$ is continuous, so $\psi=e^{i\Psi}$ with $\Psi$ continuous on $\bR$ and $\Psi(x+1)=\Psi(x)+2\pi p$, $p\in\bZ$. Then $\varphi-\Psi\circ\tilde f+\Psi$ is continuous with values in a coset of $2\pi\bZ$, hence constant; since $\Psi\circ\tilde f-\Psi$ changes by $2\pi p(\ell-1)$ under $x\mapsto x+1$ and $\varphi$ is periodic, $p=0$. Hence $\varphi$ is cohomologous to a constant.

\emph{Uniqueness of the acip and \textbf{(A6)}.} Let $g\in L^2(\bT^2)$ satisfy $g\circ F=e^{i\alpha}g$, and expand $g(x,\theta)=\sum_mg_m(x)e^{2\pi im\theta}$. Then $g_m\circ f\,e^{2\pi im\omega}=e^{i\alpha}g_m$ for each $m$. If $g_m\neq0$, then $|g_m|$ is $f$-invariant, hence constant, and by the previous step $2\pi m\omega$ is cohomologous to a constant; thus $m=0$. So $g=g_0(x)$ with $g_0\circ f=e^{i\alpha}g_0$, and since $f$ is mixing, $\alpha=0$ and $g$ is constant. Therefore Lebesgue measure on $\bT^2$ is ergodic and weakly mixing for $F$; being equivalent to Lebesgue, it is the unique acip of $F$, and \textbf{(A6)} holds.

\emph{Condition \textbf{(B2)}.} Suppose $e^{i\nu\tau}=c\,\phi\circ F/\phi$ Lebesgue-a.e., with $\nu\ne0$ and $\phi$ measurable, and expand $\phi(x,\theta)=\sum_m\phi_m(x)e^{2\pi im\theta}$. Comparing the Fourier coefficients in $\theta$ gives $e^{i\nu\tau}\phi_m=c\,e^{2\pi im\omega}\,\phi_m\circ f$ for each $m$. Choose $m$ with $\phi_m\ne0$; then $|\phi_m|$ is constant and $e^{i(\nu\tau-2\pi m\omega)}=c\,\psi\circ f/\psi$ with $|\psi|=1$. Hence $\nu\tau-2\pi m\omega$ is cohomologous to a constant, which contradicts (iii) of Lemma \ref{lem:skew} with $b=2\pi m/\nu$.

{\emph{Virtual expansion.} We also verify the VE assumption. Here $J_\frh=\ell^{-n}$ and
\[
(D_{\frh(q)}F^n)^*(a,b)=\big(\ell^n(a+s_gb),b\big).
\]
The slopes $s_g$ are uniformly bounded. By the vertical-line case of \eqref{eq:skew-nc} and \eqref{eq:skew-ss2}, for every unit covector $(a,b)$ the total weight of branches satisfying $|a+s_gb|<\ell^{-n/2}$ is at most $C\ell^{-\beta_0 n/2}$. Indeed, if $|b|$ is sufficiently small there are no such branches for large $n$; otherwise this is an interval condition on $s_g$ of width at most $C\ell^{-n/2}$. The inverse derivatives of this skew product are uniformly bounded. On the other branches, the above covector has norm at least $\ell^{n/2}$. Thus, for every fixed $s>0$,
\[
B_s(F^n)\le C_s\big(\ell^{-\beta_0 n/2}+\ell^{-sn/2}\big).
\]
The constants $\lambda_c^+>1$ and $\lambda_c^-<1$ in \textbf{(A2)} may be chosen arbitrarily close to $1$, since the central inverse expansion is bounded above and below uniformly in $n$. Choosing them so that $\lambda_c^s\ell^{-\min\{\beta_0,s\}/2}<1$ proves \eqref{virtually}. Finally, \eqref{eq:Linfty-growth} holds with $\Theta=1$, so \textbf{(P)} is satisfied. Theorem~\ref{thm:DI} and Theorems~\ref{teoA}--\ref{teoB} therefore apply; centring holds because $\int\tau\,\dd\operatorname{Leb}=0$.}
\end{ex}

\section{Useful estimates for SVPH}\label{sec:previous}

We collect some consequences of results from \cite{CaLi} for the SVPH $F$ and its transfer operator $\cL_F$. The estimates recalled below require only the domination
condition \textbf{(A4)}. The stronger pinching imposed in
\cite{CaLi} is used for higher derivative estimates, which
are not needed here.

\begin{itemize}
\item By \cite[Prop.\ 5.6]{CaLi}, there exists a constant $C_*\ge1$ such that for every $z\in \mathbb{T}^2$, any $n\in\bN$, and any vectors $v^u\in \fC_u$, $v^c\in \fC_c$ with $(a,b):=D_zF^{n}v^c\not\in \fC_u$,
\begin{equation}\label{det-form}
 C_*^{-1}\frac{\|D_zF^n v^u\|}{\|v^u\|}\frac{ |b|}{\|v^c\|}\le |\det D_zF^n|\le C_* \frac{\|D_zF^n v^u\|}{\|v^u\|}\frac{ |b|}{\|v^c\|}.
\end{equation}
In particular, for each $n\in \bN$ and $z\in \bT^2$,
\begin{equation}\label{det-est}
C_*^{-1} {\lambda_c}^{-n} \lambda^{n}\le |\det D_zF^n|\le C_* \Lambda^n {\lambda_c}^n.
\end{equation}

\item By \cite[(5.17)]{CaLi}, for each $n\in\bN$,
\begin{equation}\label{DF-1}
\|(DF^n)^{-1}\|_{\cC^0(\bT^2)}\le C\lambda_c^n.
\end{equation}

\item By \cite[Lemma 5.18 and Corollary 5.20]{CaLi},
\begin{equation}\label{C1norm}
\|\cL_F^n 1\|_{\cC^0(\bT^2)}\le C{\lambda_c}^n,\qquad
\|\cL_F^n 1\|_{\cC^1(\bT^2)}\le C{\lambda_c}^{2n}.
\end{equation}
\end{itemize}

\end{document}